\documentclass[12pt]{amsart}
\usepackage{amsmath,amsxtra,amsfonts} 
\usepackage{amssymb,amsthm,amscd,amsopn}
\usepackage{latexsym,array}
\usepackage[mathscr]{eucal}
\usepackage{graphicx,psfrag}
\usepackage{epsfig}

\newtheorem{theorem}[subsection]{Theorem}
\newtheorem{proposition}[subsection]{Proposition}
\newtheorem{lemma}[subsection]{Lemma}

\def\al{\alpha}
\def\de{\delta}
\def\eps{\varepsilon}
\def\lam{\lambda}

\def\vphi{\varphi}

\def\De{\Delta}

\def\sh{\sinh}
\def\ch{\cosh}

\def\F{{\cal F}}

\newcommand{\ud}{\mathrm{d}}

\newcommand{\cal}{\mathcal}

\newcommand{\vs}{\vspace}
\newcommand{\nd}{\noindent}
\newcommand{\hB}{\hfill$\Box$}

\newcommand{\pf}{{\bf Proof.}\, \, }
\newcommand{\rk}{{\bf Remark.}\, \, }
\newcommand{\defi}{{\bf Definition.}\ \ }

\newcommand{\ca}{\emph{Case}\;}

\newcommand{\Z}{\mathbb{Z}}
\newcommand{\R}{\mathbb{R}}

\newcommand{\N}{\mathbb{N}}

\begin{document}

\title[Besov spaces for Schr\"odinger Operators]
      {Besov spaces for Schr\"odinger Operators with barrier potentials} 
\author[J.J.~Benedetto ]
       {John J.~Benedetto}
\address{Department of Mathematics \\
        University of Maryland \\
         College Park, Maryland 20742}
 \email[J.J.~Benedetto]{jjb@math.umd.edu}      
   \urladdr[J.J.~Benedetto]{http://www.math.umd.edu/\textasciitilde{jjb}}
\author[S.~Zheng ]
       {Shijun Zheng}
\address{Department of Mathematics \\
        Louisiana State University \\
        Baton Rouge, LA 70803}
 \email[S.~Zheng]{szheng@math.lsu.edu}
 \urladdr[S.~Zheng]{http://www.math.lsu.edu/\textasciitilde{szheng}}
 \thanks{The research of the first author is supported in part by NSF/ONR.
         The second \qquad \quad author supported in part by DARPA}         
\date{\today} 

\begin{abstract}
Let $H= -\triangle +V$ be a Schr\"odinger operator on
the real line, where $V=\varepsilon^2 \chi_{[-1,1]}$. 
We define the  Besov spaces for $H$ by developing the 
associated Littlewood-Paley theory.  
This theory depends on the decay estimates of the spectral operator
$\vphi_j(H)$ in the high and low energies.
We also prove a Mikhlin-H\"ormander type multiplier theorem on these spaces, 
including the $L^p$ boundedness result.
Our approach has potential
 applications to other  Schr\"odinger operators with short-range potentials,
as well as in higher dimensions.
\end{abstract}

\maketitle 


\section{Introduction}
Let $H = - \triangle +V$ be a Schr\"{o}dinger operator on $\R$, where the potential $V$
is real-valued and belongs to $L^1\cap L^2$.
$H$ is the  Hamiltonian in the corresponding time-dependent
Schr\"{o}dinger equation 
\begin{align}
 &i \; \partial_t \psi= H \psi,\label{Eq:Schroe}\\
&\psi(0,x)=f(x) \in \mathcal{D}(H), \notag
\end{align}
where the solution is uniquely determined by the initial state: $
\psi(t,x)=e^{-i t H }f(x)$,  \mbox{$\;t\geq 0$.}   

In \cite{JN} Jensen and Nakamura introduced Besov spaces associated with
$H$ on $\R^d$ and showed that $e^{-itH}$ maps $B_p^{s+2\beta,q}(H)$ into 
$B_p^{s,q}(H)$ if $s\ge 0$, $1\le p,q\le \infty$ and 
$\beta >d\vert \frac{1}{2}-\frac{1}{p}\vert$, under the condition that
$V=V_+-V_-$ so that $V_+\in K_d^{loc}$ and
$V_-\in K_d$, $K_d$ being the Kato class.  
In this paper we generalize the definition of Besov spaces to $\al\in \R$, $0<p, q<\infty$
and show, in the case of barrier potential, that such a definition is 
independent of the choice of the dyadic system $\{\Phi,\vphi_j\}$. 

For $\alpha \in \R$, $0<p<\infty$, $0<q\leq \infty$, define the
quasi-norm  for $f\in L^2$ as
\begin{equation}\label{besov-norm}
\Vert f\Vert_{B_p^{\alpha,q}} := \Vert f\Vert_{B_p^{\alpha,q}(H)}^{\vphi} \\
= \Vert \Phi(H)\Vert_p + \left\{ \sum_{j=1}^{\infty} (2^{j\alpha} \Vert \vphi_j(H)f \Vert_p)^q
\right\}^{1/q}.
\end{equation}

\noindent
The {\em Besov spaces associated with} $H$, denoted by $B_p^{\alpha,q}(H)$
is defined to be the completion of the subspace $ L_0^2:= \{f\in L^2:
\Vert f\Vert_{B_p^{\alpha,q}} <\infty \}$ of $L^2$.

As in the Fourier case and the Hermite case \cite{Tr,Tr2,E,E2} 
we address the Besov space theory associated with $H$
by considering the Schr\"{o}dinger operator $H=-\triangle + V$, 
where  $V={\varepsilon}^2 \chi_{[-1,1]}, \varepsilon
>0$ (called the barrier potential) is one of the simplest discontinuous
potential models in quantum mechanics.  

Peetre's maximal operator plays an important role in the theory of
function spaces.  In order to establish a Peetre type maximal
inequality for $H$, we need the decay estimates of the kernel of
$\vphi_{j}(H) $ as well as of its derivative.  Based on  an
integral expression of this kernel we  obtain the decay estimates
by exploiting   the analytic behavior of the eigenfunctions
 $ e(x, \xi)$ as $\xi $ approaches $\infty $ (high energy) and $0$ (low
energy) in various cases.  When the support of $\Phi$ contains the
origin, we are in the so-called ``local energy'' case, which usually
is harder to deal with for general potentials. 
 We use certain ``matching'' method to put together integrals of the
``same type'', so that each of the resulting integrals is the Fourier transform
of a Schwartz function.   This method seems interesting and 
may have applications to other potentials.

Our first main result (Theorem \ref{th:besov}) is an equivalence theorem for
$B_p^{\alpha,q}(H)$, which tells that $\Vert f\Vert_{B_p^{\alpha,q}}^{\phi}
$ and $\Vert f\Vert_{B_p^{\alpha,q}}^{\psi}$ are equivalent
quasinorms on $B_p^{\alpha,q}(H)$, where $\{\phi_j\}, \{\psi_j\}$ are
two given smooth dyadic systems.



Using functional calculus, Jensen and Nakamura \cite{JN,JN2} obtained smooth
multiplier results for general potentials.  For barrier potential we
prove a sharp spectral multiplier theorem on $B_p^{\alpha,q}(H)$ 
(Theorem \ref{th:L^p-mult} and 
Theorem \ref{th:Besov-mult}). 



The remaining part of the paper is organized as follows. In $\S$2
we give explicit solutions to the eigenfunction equation. 
In $\S$3 we give norm characterization of $B_p^{\alpha,q}(H)$ 
using Peetre type maximal functions. The proof is
based on the decay estimates for the kernel of 
$\vphi_j(H)$.  A detailed proof of these decay estimates are included in 
$\S$4 and $\S$5.  In $\S$6 we prove a H\"ormander type multiplier
theorem for $H$. In $\S$7, we give identifications of these new
spaces $B_p^{\alpha,q}(H)$ 
with the ordinary Besov spaces for certain range of parameters $\alpha, p, q$.

\noindent
{\em Acknowledgment}\quad
The authors would like to thank A.~Jensen for his 
useful comments on the identification of Besov spaces. 








\section{Preliminaries}
2.1. \emph{Kernel formula for the spectral operator }

Let $e_+(x,\xi)$ and $e_-(x,\xi)$ be two solutions of the 
equation 
\begin{equation}\label{eq:H-eigen}
He(x,\xi) = \xi^2 e(x,\xi)
\end{equation}
with asymptotic behavior for $\xi>0$ and $\xi<0$ respectively, 
\begin{equation}\label{eq:eig-asymp}
e_\pm(x,\xi)\rightarrow \left\{
\begin{array}{ll}
T_\pm(\xi)e^{i \xi x}& x\rightarrow \pm\infty\\
e^{i \xi x} + R_\pm(\xi)e^{-i \xi x} &  x\rightarrow \mp\infty. 
\end{array}\right.
\end{equation}
Then $e_\pm (x,\xi)$ are unique for $\xi\in \R$, and equation 
(\ref{eq:H-eigen})
together with condition (\ref{eq:eig-asymp}) is equivalent to the integral equation
\begin{equation}
e(x,\xi)= e^{i\xi x} +(2i |\xi|)^{-1} \int e^{i |\xi| |x-y|} V(y)e(y,\xi)dy.
\label{LS_eq}
\end{equation}
These {\em generalized eigenfunctions} have a physical interpretation in quantum
mechanics, where $\xi^2$ is viewed as a energy parameter; they
represent the transmission and reflection waves when a particle passes
through the potential.  The coefficients $T, R$ are called the
\emph{transmission coefficient} and the {\em reflection coefficient} 
(cf. \cite{GH},  p.4179, also \cite{Fl}).
Under the condition that $V$ is in $L^1\cap L^2$, 
 we show in \cite{Zhen, Zhen2} that, 

\vspace{.122in}
\noindent
a) The essential spectrum of $H$ is $[0,\infty)$; more precisely,
$H$ has only absolutely continuous spectrum (the singular continuous
spectrum being empty); the discrete spectrum of $H$ is at most
countable. Hence if denoting $L^2$ by ${\cal H}$ we have  
${\cal H}={\cal H}_{ac}\oplus {\cal H}_{pp}$.

\vs{.104in}
\nd
b) Define the {\em generalized Fourier transform} $\mathcal{F} $ on
$L^2$: 
$$ \mathcal{F}f(\xi) : { = } \textrm{l.i.m}. \frac{1}{\sqrt{2\pi}}
\int f(x)\overline{e}(x,\xi) dx.$$
\noindent
Then  ${\cal F}$ is a unitary operator 
from ${\cal H}_{ac}$ onto $L^2$ and its adjoint
is given  by for  $g\in L^2$
\begin{displaymath}
 {\mathcal{F}}^{*} g(x) :{ = } \textrm{l.i.m}. \frac{1}{\sqrt{2\pi}}
\int g(\xi) e(x,\xi)d\xi. 
\end{displaymath}
Therefore $\vphi(H)\vert_{\cal{H}_{ac}}=\F^{*}\vphi(\xi^2)\F$.
If $H$ has no point spectrum and all generalized eigenfunctions
of $H$ are uniformly bounded in $x$ and $\xi$,
then for $f\in L^2$, 
\begin{equation}
\vphi(H) f(x) = \int  \varphi(H)(x,y) f(y)dy,
\label{ker_eq}
\end{equation}
where
\begin{equation}\label{eq:K(x,y)}
\vphi(H)(x,y)=\frac{1}{2\pi} \int_{-\infty}^{\infty} \vphi({\xi}^2) 
 e(x,\xi)\overline{ e}(y,\xi) d \xi.
\end{equation}




A variant of the formula (\ref{ker_eq}) can be found in \cite{GH} 
for short-range potentials defined as a measure. In 3D similar
formula is used by Tao \cite{Tao} in a scattering problem.

Since the barrier potential $V=\eps^2 \chi_{[-1,1]}$ is in $L^1\cap L^2$ 
and the eigenfunctions 
of $H$ are uniformly bounded (see subsection 2.3), the formula
(\ref{ker_eq}) is valid for $V$. Note that the
corresponding point spectrum is empty.
 

2.2. \emph{Dyadic system and Besov spaces} 

Let $\Phi, \vphi, \Psi, \psi $ be $C^{\infty}$ smooth functions, satisfying
 
$i$) $\textrm{supp}\; \Phi,\; \textrm{supp}\; \Psi \subset\{ |\xi|\leq 1\}; \;
|\Phi (\xi)|,\; |\Psi (\xi)| \geq c>0 $ if  $|\xi|\leq \frac{1}{2}$;

$ii$) $\textrm{supp} \; \vphi,\; \textrm{supp}\; \psi \subset \{\frac{1}{4}\leq |\xi|\leq
1 \};$

$|\vphi ( \xi ) |,\; |\psi ( \xi )|\geq c>0$  if $\frac{3}{8} \leq |\xi|\leq
 \frac{7}{8};$

$iii$) $\; \Phi (\xi) \Psi (\xi) + \sum_{j=1}^{\infty} \vphi(2^{-j} \xi
)\psi (2^{-j} \xi)
=1,  \;\;\forall \;\xi \in \R $. 

\nd
The existence of 
such functions can be justified by consulting e.g., \cite{FJW}.
The almost orthogonal relation $(iii)$ for the system 
$\varphi_j(\xi):=\varphi(2^{-j}\xi)$,
$\psi_j(\xi):=\psi(2^{-j}\xi)$  allows us to give an expansion of $f$ in $L^2$ as 
$$
f= \Phi(H) \Psi(H)f + \sum_j \vphi_j(H) \psi_j(H)f.$$

As in the Fourier case, let $0 < p,q \le \infty, \alpha \in
\R$ we define the $B_p^{\alpha,q}$ quasi-norm as in (\ref{besov-norm}) for
$f\in L^2$.
Note that when $0<p<1$ or $0<q<1$, we can always define a metric $d$ on
$B_p^{\alpha,q}$, so that the metric space $( B_p^{\alpha,q}, d) $ is
topologically isomorphic to the quasi-normed space.
In fact, Lemma 3.10.1 in \cite{BL} tells that every quasi-normed
linear space is metrizable.


  
2.3. \emph{Generalized eigenfunctions of $H$}

We now determine 
eigenfunctions of $H=-\triangle +V$, where $V=
\varepsilon^2 \chi_{[-1,1]} $, also see, e.g., \cite{Fl}. 

First $e(x,\xi)$ must have the following form.
If  $\xi >\varepsilon$, then
$$
 e(x,\xi)= \left\{
\begin{array}{ll}
A_{+}e^{i \xi x} +A^{\prime}_{+}e^{-i \xi x}& x<-1\\
B_{+}e^{i K x} +B'_{+}e^{-i K x}& |x|\leq 1\\
C_{+}e^{i \xi x} +C'_{+}e^{-i \xi x}& x>1,
\end{array}\right.
$$
where $K=\sqrt{ {\xi}^2 -{\varepsilon}^2 }$; and if $0<\xi <\eps$,
then 
$$
 e(x,\xi)= \left\{
\begin{array}{ll}
A_{+}e^{i \xi x} +A'_{+}e^{-i \xi x}& x<-1\\
B_{+}e^{\rho x} +B'_{+}e^{-\rho x}& |x|\leq 1\\
C_{+}e^{i \xi x} +C'_{+}e^{-i \xi x}& x>1,
\end{array}\right.
$$
where $\rho=\sqrt{{\eps}^2 -\xi^2}$. 

The Lippman-Schwinger equation (\ref{LS_eq}) requires that  
$ e(x,\xi)$ is differentiable  in $x$. By the $C^1$ condition at $\pm
1$ we can obtain the
precise values of the coefficients $A, A', B, B', C,C'$ as follows.

Let 
$$ \rho =\rho(\xi)= \left\{ 
\begin{array}{ll}
i K=i\sqrt{\xi^2 -{\varepsilon}^2} & |\xi|>\varepsilon \\
\sqrt{{\varepsilon}^2 -\xi^2}&   |\xi|\leq \varepsilon.
\end{array}\right.
$$

Then for $\xi>0$,
\begin{align*}
C'_{+}=&\; 0,  \;\;A_{+}=1,\\
C_{+}=&\frac{2\rho \xi e^{-2 i\xi}}{2\rho \xi \cosh 2\rho + i(\rho^2-\xi^2)
\sinh 2\rho}\\
A'_{+} =&-i \frac{C_+}{2\rho \xi}\varepsilon^2 \sinh
2\rho=-i\frac{\eps^2 \sinh 2\rho  e^{-2 i\xi}}{2\rho \xi \cosh
2\rho + i(\rho^2-\xi^2)\sinh 2\rho}\\
B_+=&\;\frac{C_+}{2\rho}(\rho+i\xi)e^{-\rho+i\xi}, \qquad \;\;
B'_{+}=\frac{C_+}{2\rho}(\rho-i\xi)e^{\rho+i\xi}.
\end{align*}

For $\xi<0$, we obtain similarly, with the same notation 
$\rho=\rho (\xi)$,

$$
 e(x,\xi)= \left\{
\begin{array}{ll}
A_{-}e^{i \xi x} +A'_{-}e^{-i \xi x}& x<-1\\
B_{-}e^{\rho x} +B'_{-}e^{-\rho x}& |x|\leq 1\\
C_{-}e^{i \xi x} +C'_{-}e^{-i \xi x}& x>1,
\end{array}\right.
$$
where $C_{-}=1, \;A'_- =0$,

\begin{align*}
A_{-}=&\frac{2\rho \xi e^{2 i\xi}}{2\rho \xi \cosh 2\rho -i(\rho^2-\xi^2)
\sinh 2\rho},\\
C'_{-} =& i \frac{A_-}{2\rho \xi}\varepsilon^2 \sinh
2\rho= i\frac{\varepsilon^2 \sinh 2\rho  e^{2 i\xi}}{2\rho \xi \cosh
2\rho - i(\rho^2-\xi^2)\sinh 2\rho}\\
B_- =&\frac{A_-}{2\rho}(\rho+i\xi)e^{\rho-i\xi}, \;\;\;
{B'}_{-}=\frac{A_-}{2\rho}(\rho-i\xi)e^{-\rho-i\xi}.
\end{align*}


Furthermore, if we define for $\xi \in \R\setminus \{0\}$

$$
 e(x,\xi)= \left\{
\begin{array}{ll}
e_{+}( x, \xi) & \xi>0\\
e_{-}( x, \xi) & \xi<0,
\end{array}\right.
$$ 
then 
\begin{equation}
e( x, -\xi) = e( -x, \xi), \;\;\xi \neq 0, 
\label{eig_eq}
\end{equation}
which follows from the following simple relations between the coefficients:

$$A_{-}(\xi)=\overline{A}_{-}(-\xi ) =\overline{C}_+(\xi) =C_+(-\xi)$$

$$C'_-(\xi) =\overline{C'}_-(-\xi)= \overline{A'}_+(\xi)= A'_+(-\xi )$$
and  
$$ B_{+}(-\xi ) = B'_{-}(\xi ),\;\; B_{-}(-\xi ) = B'_{+}(\xi).$$






\vspace{0.2in}\noindent
Identity (\ref{eig_eq}) allows us to simplify the estimation in various cases; see \S 4--6.
Some of the above relations 
can also be found in \cite{GH} for general potentials.

\section{Peetre type maximal inequality}

Let $\Phi, \vphi, \Psi, \psi $ be $C^{\infty}$ smooth functions,
satisfying the conditions given in $\S$2.  Recall that if $\phi\in C_0 (\R)$, 
the  operator $\phi (H)$ 
has the kernel (\ref{eq:K(x,y)}).
Note that $e (\cdot , \xi)\in C^1 (\xi \neq 0, \pm \epsilon )$ implies
$\phi (H) (x,y)\in C^1(\R \times \R)$.
 
\begin{lemma}  Let $K_j(x,y)=\vphi( 2^{-j} H ) (x,y)$.  
$\;\textrm{supp}\; \vphi \subset [\frac{1}{4},1]$.\\
a) If $j>4 + 2\log_2 \epsilon$, 
we have for each $n\in \Z^+$,  
$$ \vert K_j(x,y) \vert \leq C_n \sum_{\ell=0}^{2N} 2^{j/2}
(1+2^{j/2} |x \pm y \pm 2\ell|)^{-n},$$
where $N=$ the smallest integer $\geq \max \{1, n/4\}$. \\
b) If $-\infty<j\leq J:=4+[2\log_2 \epsilon ] $, then for each $n\geq 0$ 
$$\vert K_j(x,y) \vert \leq C_n  2^{j/2}
(1+2^{j/2}|x \pm y|)^{-n}. \; \;$$
\end{lemma}

\begin{lemma}  Let $K(x,y)=\Phi ( H)(x,y)$, 
$supp \;\Phi\subset [-1,1]$. We have for each $n\geq 0$
$$ 
\vert K(x,y) \vert \leq C_{n}  
(1+|x - y|)^{-n}.
$$
\end{lemma}

\vspace{0.2in}
We also need decay estimates for the derivative of the kernel.

\begin{lemma} Let $\vphi_j(x)=\vphi(2^{-j}x)$. $K_j(x,y)=\vphi_j(H)(x,y)$. 
\\ 
a) If $j>J$, 
then for each $n$ there is a constant $C_n$: 
$$ \bigg{\vert} \frac{\partial}{\partial x} K_j(x,y) \bigg{\vert} \leq C_n \sum_{
\ell = 0}^{2N} 2^j
(1+2^{j/2}|x \pm y \pm 2\ell |)^{-n},$$
where $ N$ means the the same as in Lemma 3.1 (a).\\ 
b) If $-\infty<j\leq J$, then for each $n$  there is a
constant $C_n$: 
$$ \left\vert \frac{\partial}{\partial x} K_j(x,y) \right\vert \leq C_n  2^j
(1+2^{j/2}|x \pm y |)^{-n}.$$
\end{lemma}

\begin{lemma}  Let $\Phi$ be as in Lemma 3.2. Then for each $n$
$$ 
\left\vert  \frac{\partial}{\partial x} K(x,y) \right\vert \leq C_{n}  
(1+|x - y|)^{-n}.
$$
\end{lemma}

\vspace{0.2in}
Proofs of Lemma 3.1$-$3.4 are given in $\S$4 and $\S$5,
 which are elementary calculus but quite lengthy. 
These lemmas are essential for us to establish a Peetre type maximal
inequality.

Given $s>0$ define the \emph{Peetre maximal functions} for $H$ by: if $j>J$,
$$
\vphi_j^*f(x) = \sup_{t\in \R} \frac{| \vphi_j(H)f(t)|}{\min_{\ell, \pm}
 ( 1+2^{j/2}|x \pm t \pm 2\ell |)^s} , $$
and 
\begin{align*}
\vphi_j^{**} f (x) = \sup_{t\in \R} \frac{| (\vphi_j(H)f)'(t)|}{\min_{\ell, \pm}
 ( 1+2^{j/2}|x \pm t \pm 2\ell |)^s},
\end{align*}
where the minimum is taken over $0\le \ell \le 2N$ and 
$N=\max\{1,\frac{[s]+2}{4}\}$.

Similarly, define for $j\leq J$,
\begin{align*}
\vphi_j^*f(x) =& \sup_{t\in \R} \frac{| \vphi_j(H)f(t)|}{\min_{ \pm}
 ( 1+2^{j/2}|x \pm t  |)^s} , \;\;s>0,\\
\Phi^*f(x) =& \sup_{t\in \R} \frac{| \Phi(H)f(t)|}{
 ( 1+|x - t  |)^s},\\
\vphi_j^{**}f(x) =& \sup_{t\in \R} \frac{| (\vphi_j(H)f)'(t)|}{\min_{ \pm}
 ( 1+2^{j/2}|x \pm t  |)^s},\\
\Phi^{**}f(x) =& \sup_{t\in \R} \frac{| (\Phi(H)f)'(t)|}{
 ( 1+|x - t  |)^s}.
\end{align*}
 We have used the abbreviation $\vphi_j^*f:=\vphi_{j,s}^*f$.
Notice that 
\begin{equation}\label{ineq:phi*-phi}
\vphi_j^*f(x) \geq |\vphi_j(H)f(x)|.
\end{equation}

In the following we slightly abuse the notation 
 $\vphi_0^*f=\Phi^*f$, etc, in case of no confusion.
\begin{lemma} For $s>0$, there exists a constant $c_s>0$ such  that
$$
\vphi_j^{**}f(x) \leq c_s 2^{j/2} \max_{k,\pm}
 \vphi_j^{*} f( x\pm 2k),$$
where the maximum is taken over $0\leq k \leq 2N$ and both $\pm$. 
\end{lemma}

\nd
\pf From the identity
$$
\vphi_j(H)f(x) =\sum_{\nu = -1}^1  (\vphi \psi)_{j+\nu}(H) \vphi_j(H)f ,
\;\; f\in L^2$$
with convention $\vphi_0=\Phi$ and $\vphi_{-1}=0,$
we derive 
$$
\frac{d}{dt}( \vphi_j(H)f)(t) 
=  \sum_{\nu = -1}^1 \int_{\R} \frac{\partial}{\partial t}K_{j+\nu}(t,y)
\vphi_j(H)f(y)dy,$$ 
where $K_j$ denotes the kernel of $(\vphi\psi)(2^{-j}H)$. 

Let $j>J$ first.  Apply Lemma 3.3  to get 
$$
\frac{ |\frac{d}{dt}( \vphi_j(H)f)(t)|}{\min_{k, \pm} 
( 1+2^{j/2}|x \pm t \pm 2k |)^s} \leq
C_n 
 \sum_{\nu=-1}^1  \sum_{\ell, \sigma,\mu} \int_\R $$

$$
\frac{2^{j+\nu}}{(
1+2^{\frac{j+\nu}{2}} |t +\sigma y +\mu 2\ell |)^n} \times
\frac{ |\vphi_j(H)f(y)| dy }{\min_{k, \pm} 
( 1+2^{j/2}|x \pm t \pm 2k |)^s} ,$$
where the inner sum is taken over all $0\le \ell\le 2N$ and $\sigma, \mu \in \{\pm 1\}$; similar notation for $\min_{\ell, \pm}$.

\vspace{0.15in}

\nd
\emph{Claim}. 
$$\frac{ |\vphi_j(H)f(y)| }{\min_{\ell, \pm} 
( 1+2^{j/2}|x \pm t \pm 2\ell |)^s} 
\leq \max\limits_{k,\pm} \vphi_j^*f(x \pm 2k) \min \limits_{\ell,\pm}
( 1+2^{j/2}|t \pm y \pm 2\ell |)^s .$$ 

To prove the claim, note that $\exists\; \de,\; \eps\in\{\pm 1\}  $ and $\ell_0$
such that   $\min_{\ell,\pm}
( 1+2^{j/2}|x \pm t \pm 2\ell |) =  1+2^{j/2}|x + \delta t +\epsilon
2\ell_0 | $ 
for given $x$, $t$. Then for each $\sigma$, $\mu$ and $\ell$ the left hand
side is less than or equal to 
\begin{align*} 
&\frac{ |\vphi_j(H)f(y)| }{\min_{k, \pm} 
( 1+2^{j/2}|x +\epsilon \cdot  2\ell_0 \pm y \pm 2k |)^s} \cdot
\frac{  (1+2^{j/2}|x +\epsilon \cdot  2\ell_0 +\sigma' y +
\mu' 2\ell |)^s}
{( 1+2^{j/2} |x +\delta t  +\epsilon 2\ell_0 |)^s}\\
&\leq 
\vphi_j^*f(x +\epsilon 2\ell_0) 
(1+2^{j/2}|-\delta t +\sigma' y +
\mu' 2\ell |)^s \\
&\le  \max_{k, \pm}\vphi_j^*f(x \pm 2k)
 (1+2^{j/2}|t +\sigma y +\mu 2\ell |)^s,
\end{align*}
where we put ${\sigma}'= -\delta \sigma, {\mu }'= -\delta \mu$ and used
for $s>0$, 
$$(1+2^{j/2} | x  + \epsilon 2\ell_0 +\sigma' y +\mu' 2\ell |)^s
\leq ( 1+2^{j/2} |x +\delta t + \epsilon 2\ell_0 |)^s (1+2^{j/2}|-\delta t
+\sigma' y +\mu' 2\ell |)^s .$$
Since $\sigma, \mu$, $\ell$  are arbitrary, 
the claim is proved.

It follows that 
\begin{align*}
 &\frac{ |\frac{d}{dt}( \vphi_j(H)f)(t)|}{\min_{k, \pm} 
( 1+2^{j/2}|x \pm t \pm 2k |)^s} \leq C_n
  \sum_{\nu=-1}^1 \sum\limits_{\ell, \sigma,\mu } 
\max_{k, \pm}\vphi_j^*f(x \pm 2k) \times \\
&\qquad\qquad\int_{\R} \frac{2^{j+\nu}}{(
1+2^{\frac{j+\nu}{2} } |t +\sigma y +\mu 2\ell |)^n }
\cdot 
( 1+2^{j/2}|t +  \sigma y +\mu 2\ell |)^s \ud y \\
\leq& C_{n} \max_{\substack{0\leq k\leq 2N\\* \pm}} \vphi_j^*f(x \pm
2k)\sum\limits^{2N}_{\substack{\ell =0\\* \sigma,\mu=\pm 1}} 
\int_{\R} 
\frac{2^{j+n/2}}{(1+2^{j/2} |t +\sigma y +\mu 2\ell |)^{n-s}}
dy \\
\leq &C_{n,s} (2N+1)\max_{\substack{0\leq k\leq 2N\\ \pm}}
\vphi_j^*f(x \pm 2k) 2^{j/2},
\end{align*}
provided $n-s>1$. 
Thus one may take $n=[s]+2$. 

For $j\le J$ similarly we obtain the following inequalities, using Lemma
3.3(b) and Lemma 3.4 in place of Lemma 3.3(a), 
$$ \vphi_j^{**} f(x) \le C 2^{j/2} \vphi_j^{*} f(x)$$
and 
$$ \Phi^{**} f(x) \le C  \Phi^{*} f(x).$$ 
This proves Lemma 3.5. \hfill$\Box$

We are ready to show Peetre maximal inequality for $H$. Let 
$M$ be the Hardy-Littlewood maximal operator:
$$
Mf(x)=\sup_{x\in I} |I|^{-1}\int_{I}|f(u)|du,   
$$
where the supreme is taken over all intervals $I$ containing $x$.

\begin{lemma} Let $0<r<\infty$. There exists a constant $C>0$
independent of $0<\varepsilon \leq 1$ such that
$$
 \vphi_j^*f(x) \leq C \epsilon \sum_{\ell =0}^{2N}\vphi_j^*f(x \pm 2\ell)
+C {\epsilon}^{-1/r} \sum_{\ell =0}^{2N} [M({|\vphi_j(H)f|}^r)]^{1/r}(\pm x
\pm 2\ell),$$
where $\epsilon >0$ can be chosen arbitrarily small.
\end{lemma}

\nd
Remark 1. It is well known that
$M$ is bounded on $L^p, 1<p<\infty$.  Lemma 3.6 implies that if $s=1/r$, then

\begin{equation}\label{eq:phi*-Lp}
\Vert  \vphi_j^*f \Vert_p \leq c \Vert  \vphi_j(H)f \Vert_p, \;\; 0<p \le \infty
\end{equation}
by taking $\epsilon$ small enough and  $0<r <p \;(s=1/r >1/p)$.

\nd
Remark 2. For $j\le J$, the inequality in Lemma 3.6 takes a simpler form 
\begin{align*}
\vphi_j^* f(x) \le& C_s \epsilon^{-s} [M(|\vphi_j(H)f |^r )]^{1/r}(\pm x),\\  \Phi^* f(x) \le& C_s \epsilon^{-s} [M(|\Phi(H)f |^r )]^{1/r}( x).
\end{align*}
Compare the analogue in Fourier case \cite{Tr} and Hermite case 
\cite{E}.

\vspace{.15in}
\nd
\pf Let $g(x) \in C^1(\R)$. As in \cite{Tr}, the mean value theorem
gives for $z_0\in \R, \;\delta>0$ 

$$
|g(z_0)|\leq 2\delta \sup_{\substack{|z-z_0|\leq \delta}} |g'(z)|
 +(2\delta)^{-1/r}\left( \int_{|z-z_0|\leq \delta} |g|^r dz \right)^{1/r}.$$

\nd
Put $g(z)=\vphi_j(H)f(\pm x \pm 2\ell -z)\in C^1$ to get, with
$0<\delta\leq 1, 0\leq \ell \leq 2N$,
\begin{align*}
&\frac{|\vphi_j(H) f (\pm x \pm 2\ell -z)|}{(1+2^{j/2}|z|)^{1/r}}
\leq 2\delta \sup_{\substack{|u-z|\leq \delta}} \frac{(1+2^{j/2}|u|)^{1/r}
|\frac{d}{dz}(\vphi_j(H)f)(\pm x \pm 2\ell -u)|}{(1+2^{j/2}|z|)^{1/r}
(1+2^{j/2}|u|)^{1/r}} \\
&+ (2\delta)^{-1/r} (1+ 2^{j/2}|z|)^{-1/r} \left(\int_{|u-z|\leq \delta} 
|\vphi_j(H)f(\pm x \pm 2\ell -u)|^r du\right)^{1/r}\\
&\leq 2\delta (1+2^{j/2}\delta)^{1/r} \sup_{\substack{u\in \R}}
\frac{ | \frac{d}{dz}(\vphi_j(H)f)(\pm x \pm 2\ell -u)|}{(1+2^{j/2}
|u|)^{1/r}}\\
&+ 
(2\delta)^{-1/r} (1+2^{j/2}|z|)^{-1/r} \left(\int_{|u|\leq |z|+\delta} 
|\vphi_j(H)f(\pm x \pm 2\ell -u)|^r du\right)^{1/r}\\
&= 2\delta (1+2^{j/2}\delta)^{1/r} \vphi_j^{**}f(x) +{\delta}^{-1/r} 
\left(\frac{|z|+\delta}{1+2^{j/2}|z|}\right)^{1/r} [M(|\vphi_j(H)f(\pm x \pm 2\ell
)|^r )]^{1/r}\\
&=C \epsilon \sum_{\ell =0}^{2N} \vphi_j^*f(x\pm 2\ell) +(1+ \epsilon^{-1})^{1/r}
\sum_{\ell =0}^{2N} [M(\vphi_j(H)f)^r]^{1/r} (\pm x \pm 2\ell),
\end{align*}
by taking $\delta= 2^{-j/2}\epsilon $ and using Lemma 3.5. This
proves the lemma. \hfil$\Box$ 

A direct consequence of Lemma 3.6 is the Peetre maximal function
characterization of  the spaces $ B_p^{\alpha,q}(H)$.

\begin{theorem} \label{th:besov} Let $\alpha\in \R,  0<p,q \leq \infty$. 
If $\vphi_j^*f$ and $\Phi^*f$
are defined 
with $s>1/p$, we have for $f\in L^2$
\begin{equation}
 \Vert f\Vert_{B_p^{\alpha,q}} \approx \Vert {\Phi}^*f \Vert_p+
\left(\sum_{j=1}^{\infty} 2^{j\alpha q}\Vert \vphi_j^*f
\Vert_p^q \right)^{1/q}.
\end{equation}
Furthermore, $B_p^{\alpha,q}$ is a quasi-Banach space (Banach space if
$p\geq 1, q \geq 1$) and it is independent of the choice of $\{\Phi,
{\vphi}_j\}_{j\geq 1}$.
\end{theorem}

\nd
\pf In view of (\ref{ineq:phi*-phi}), it is sufficient to show that
$$ \Vert \Phi^*f\Vert_p + \left(
\sum_{j=1}^{\infty}2^{j \alpha q} \Vert \vphi_j^* f \Vert_p^q\right)^{1/q} \leq C
\Vert f \Vert_{ B_p^{ \alpha , q } }
$$
but this follows from (\ref{eq:phi*-Lp}) immediately.

Next we show that $B_p^{\alpha,q}$ is independent of the generating
functions, i.e., 
given two systems $\{\phi_j,\psi_j \}$ and 
$\{\tilde{\phi}_j,\tilde{\psi}_j\} $, then
$\Vert f\Vert^{\phi}_{B_p^{\alpha,q}}$ and $\Vert
f\Vert^{\tilde{\phi}}_{B_p^{\alpha , q } }$ are equivalent quasi-norms on $B_p^{\alpha,q}$.

Write $\phi_j(H)= \sum^1_{\nu =- 1}
\phi_j(H)(\tilde{\phi}\tilde{\psi})_{j+\nu}(H)$
by the identity $\phi_j(x)=\phi_j(x) \sum^1_{\nu =-1}
(\tilde{\phi}\tilde{\psi})_{j+\nu}(x)$, 
$\forall \;x$. 
We have by Lemma 3.1,
\begin{align*}
& \vert \phi_j(H)f(x) \vert \leq \sum^1_{\nu =- 1}
 \sum_{\ell,\pm}\int_\R \frac{2^{j/2}}{ (1+2^{j/2} |x \pm y \pm
2\ell |)^n} \vert \tilde{\phi}_{j+\nu}(H)f(y) \vert dy\\
\leq& C\sum^1_{\nu =- 1}  \tilde{\phi}_{j+\nu}^*f(x) 
\sum_{\ell,\pm}\int_\R
\frac{2^{j/2}}{ (1+2^{j/2} |x \pm y \pm 2\ell |)^n}
\min_{\substack{k, \pm }} 
(1+2^{j/2} |x \pm y \pm 2k |)^s dy\\
\le&  C_N \sum^1_{\nu =- 1}  \tilde{\phi}_{j+\nu}^* f(x),
\end{align*}
provided $n-s>1$. Thus for $f\in L^2$
$$
\Vert f\Vert^{\phi}_{B_p^{\alpha,q}}= \Vert \{2^{j\alpha}\Vert
\phi_j(H)f\Vert_p\} \Vert_{\ell^q} \leq C_{\alpha} \Vert \{2^{j\alpha}\Vert
\tilde{\phi}_j^*f\Vert_p\} \Vert_{\ell^q} \approx  \Vert
f\Vert^{\tilde{\phi}}_{B_p^{\alpha,q}}.$$
This concludes the proof of Theorem 3.7 (That $B_p^{\alpha,q}$ are
quasi-Banach spaces follows directly from the definition).   $\hfill \Box$

As expected from Lemma 3.6 we can define the homogeneous Besov spaces 
and obtain a maximal function
characterization as well.

Let ${\vphi,\; \psi}\in C^{\infty}$  satisfy
 
$i$) $\textrm{supp} \; \vphi,\; \textrm{supp}\; \psi \subset \{\frac{1}{4}\leq |\xi|\leq
1 \};$

$|\vphi ( \xi ) |, \; |\psi ( \xi )|\geq c>0$  if $\frac{3}{8} \leq |\xi|\leq
 \frac{7}{8};$

$ii$) $\;  \sum_{j=-\infty}^{\infty} \vphi(2^{-j} \xi
)\psi (2^{-j} \xi)
=1,  \;\;\forall \;\xi \neq 0 $. 

\nd
\defi The {\em homogeneous Besov space}
${\dot{B}}_p^{\alpha,q}:= {\dot{B}}_p^{ \alpha,q } (H)
$ associated with $H$ is the completion of the set $\{f\in L^2: 
\Vert f\Vert_{ {\dot{B}}_p^{ \alpha,q } } <\infty\}$ with respect to the
norm $\Vert \cdot \Vert_{{ \dot{B}}_p^{\alpha,q}}$,
where 
$$\Vert f\Vert_{ {\dot{B}}_p^{\alpha,q} } = \left\{\sum_{j=-\infty}^{\infty}
(2^{j\alpha}\Vert\vphi_j(H)f \Vert_p)^q \right\}^{1/q}
.$$

\begin{theorem} Let $\alpha\in \R,  0<p,q \leq \infty$. If $\vphi_j^*f$
is defined for $j\in \Z $ with $ s>1/p $, then for $f\in L^2$
$$
 \Vert f\Vert_{ {\dot {B} }_p^{ \alpha , q }} \approx \left( 
\sum_{j=-\infty}^{\infty}  2^{j\alpha q}\Vert \vphi_j^*f
\Vert_p^q \right)^{1/q}.$$
Furthermore, $\Vert \cdot \Vert^{\phi}_{\dot{B}_p^{\alpha,q}}$ and
$\Vert \cdot \Vert^{\tilde{\phi}}_{\dot{B}_p^{\alpha,q}}$ are  equivalent norms
on the  quasi-Banach space $\dot{B}_p^{\alpha,q } $ for any given two
systems $\{\phi_j\}$ and $\{\tilde{\phi}_j\}$.
\end{theorem}
The proof is completely implicit in that of Theorem 3.7 and hence
omitted.

Moreover, like in the Fourier case and Hermite case \cite{Tr}, \cite{E}, 
Peetre maximal inequality enables us to define and characterize
Triebel-Lizorkin spaces; see \cite{Zhen2}. 

\section{High and low energy estimates}
We give proofs of Lemma 3.1 and Lemma 3.3 for the decay estimates of
the kernel $\vphi_j(H)(x,y)$ and  $\frac{\partial}{\partial x} \vphi_j(H)(x,y)$.
Recall that  if $H=\int \lambda dE_{\lambda}$ is the spectral resolution
of $H$, then $\vphi_j(H)=\int
\vphi (2^{-j}\lambda) dE_{\lambda}=
{\cal F}^{-1} \vphi_j (\xi^2)\mathcal{F}$   with $\textrm{supp} \;\vphi_j \subset
[2^{j-2}, 2^j ] $,
which means that the spectrum of $\vphi_j(H)$ is bounded
away from $0$.

When $j>J=4+[2\log_2 \varepsilon]$, we treat $K_j(x,y)$, the kernel
of the operator $\vphi_j (H)$,
as an oscillatory integral as $\xi\rightarrow \infty$.
When $j\leq J$, we use the asymptotic property (as $\xi \rightarrow 0$) of eigenfunctions
$e(x,\xi ) $ to get estimates for the kernel.

Since $e(x,\xi ) $ has different expressions as $x>1, \;|x|\leq 1$ and
$x<-1 $, the estimates are divided into nine cases, namely, 
\begin{align*}
&1a.\; x>1, y>1; \; 1b.\; x>1, |y|\leq 1; \;1c.\; x>1, \, y <-1 \\
& 2a. \; |x|\leq 1, y>1; \; 2b.\; |x|\leq 1, |y|\leq 1; \;2c.\; |x|\leq 1,\,y <-1\\
&3a.\;  |x|<-1, y>1; \; 3b.\; |x|<- 1, |y|\leq 1;\; 3c.\; x<- 1, y <-1.
\end{align*}

\nd
By virtue of the relation $ e (x,-\xi ) = e (-x,\xi ) $ and the trivial
conjugation relation 
$\vphi(\lambda^2 H)(x,y)  = \overline{\overline{\vphi}(\lambda^2 H)(y,x)}
=\vphi (\lam H)(-x,-y).$
we see, however, these cases reduces to the following four cases: 
$1a, 1b, 1c, 2b$. 

Let $\lambda=2^{-j/2}$, then ${\lambda}^{-1}>4\epsilon \Leftrightarrow 
j>J=[2\log_2 \varepsilon]+4$.
Recall from (\ref{eq:K(x,y)}) that
\begin{equation}\label{eq:Kj}
K_j(x,y)= \frac{1}{2\pi} \int \psi (\lambda
\xi ) e(x, \xi)\overline{e}(y,\xi ) d\xi,
\end{equation}
where $\psi (x) = \vphi (x^2)$, with $\vphi$ satisfying 
$ \vphi \in C^{\infty}_0, \textrm{supp} \;\vphi\subset 
[\frac{1}{4}, 1]\cup [-1,-\frac{1}{4}]$.

4.1. \emph{High energy estimates} $j>J$

 {\bf Proof of Lemma 3.1(a).\quad}  We only show Cases 1a and 2b.  Cases
1b and 1c can be shown similarly.  

\nd
$Case\; 1a$. $ x>1, y>1$. Let $I(x,y)=2\pi K_j(x,y)$. Then by 
(\ref{eq:Kj}) 
\begin{align*}
I(x,y) & = \int_{1/{ 2\lambda} }^{ 1/\lambda} \psi (\lambda \xi )C_+ e^{ix
\xi}\overline{C_+ e^{iy \xi}} d\xi \\
\;& + \int^{-1/{ 2\lambda} }_{- 1/\lambda}  \psi (\lambda \xi )(
e^{i x \xi} 
+{C'}_- e^{-i x \xi})\overline{ e^{iy\xi}+{C'}_- e^{-iy\xi} } d\xi
:=I^+ +I^-.
\end{align*}

Convention. $\int^+ =\int_{1/{ 2\lambda } }^{ 1/\lambda}$, \quad
 $\int^- =\int_{-1/\lambda }^{- 1/{ 2 \lambda } } $.

We break the estimate of $I^+$ into two parts:

$$ \int^+ =\int_{1/{ 2\lambda} }^{ 1/\lambda}  \psi (\lambda \xi ) |C_+|^2 
 e^{i(x-y)\xi} d\xi$$

$$
=\int_{1/{ 2\lambda}}^{ 1/\lambda} \psi (\lambda \xi )\frac{4K^2 {\xi}^2}
{4K^2 {\xi}^2+\eps^4 \sin^2 2K } e^{ i (x-y) \xi } d \xi \;\;(K=\sqrt{\xi^2-\eps^2})
$$

$$
\leq \sum_{p=0}^{N-1} \bigg{\vert} \int_{1/{ 2\lambda} }^{ 1/\lambda} \psi (\lambda \xi)
(\frac{\eps^4 {\sin}^2 2K}{4K^2 {\xi}^2} )^p  e^{i(x-y)\xi} d\xi \bigg{\vert}$$

$$
+\bigg{\vert} \int_{1/{ 2\lambda} }^{ 1/\lambda} \psi (\lambda \xi)
\tilde{O}(\frac{\eps^4 \sin^2 2K}{4K^2 {\xi}^2} )^N  e^{i(x-y)\xi}
d\xi \bigg{\vert} :=I_N^+ + R^+_N,$$
where we used
\begin{equation}
\frac{4K^2 {\xi}^2}{4K^2 {\xi}^2+\eps^4 \sin^2 2K} = \sum_{p=0}^{\infty}
(-1)^p \left( \frac{\eps^4 \sin^2 2K}{4K^2 {\xi}^2} \right)^p
\end{equation}
because $\frac{\eps^4 \sin^2 2K}{4K^2 {\xi}^2} \leq \frac{\eps^4 }{
3{\xi}^4} \leq \frac{1}{3} (\frac{1}{2})^4 <1$, if $|\xi|\geq
1/2\lambda >2\eps \;( K^2\geq  \frac{3}{4} {\xi}^2)$. 

\nd
Notation.$\; \widetilde{O}(\xi^{-m}):= \widetilde{ O_{ \infty } } ( \xi^{-m}) $
denotes a function whose derivatives of arbitrary order $\geq 0$ has
estimates $ {O}(\xi^{-m})$, as $\xi \rightarrow \infty$. 

Note that
$$\sum_{p=N}^{\infty}(-1)^p  \left( \frac{\eps^4 \sin^2 2K}{4K^2 {\xi}^2}\right)^p
=(-1)^N  \left(\frac{\eps^4 \sin^2 2K}{4K^2 \xi^2}\right)^N \bigg{/} \left(1+ \frac{\eps^4 \sin^2 2K}{4K^2 \xi^2}\right).$$

If we write $\sin 2K =(2i)^{-1} (e^{i2K}-e^{-i2K})$, the integral in
each term of the sum $I^+_N$ is bounded by a linear combination of the
form, with $0\leq \ell \leq 2p, 0\leq p \leq N-1,$ 

$${\int^+}'= \int^+ \psi(\lambda \xi)  \frac{ e^{\pm i2K \ell} }{(4K^2 \xi^2)^p}
 e^{i(x-y)\xi} d\xi$$

$$= \int^+ \psi (\lambda \xi) (4K^2\xi^2)^{-p} e^{\pm i2(K-\xi)\ell}
e^{i(x-y\pm 2\ell)\xi} d\xi.$$

The following estimates will be used often. 
$$
\left\{ 
\begin{array}{lcl}
\frac{d^n}{d \xi^n } [\psi (\lambda \xi)] &= & {\lambda}^n \psi^{(n)}(\lambda
\xi )  \leq C {\lambda}^n\\
\frac{d^i}{ d \xi^i }[ (K^2 \xi^2)^{-p}] &=&
\frac{d^i}{ d \xi^i }(\frac{1}{ \xi^{4p} [1-\eps^2/{ \xi^2} ]^p })\\
\;&=& \quad \frac{d^i}{d\xi^i } \bigg(  \sum_{k=0}^{\infty}{-p \choose k}\frac{(-1)^k
\eps^{2k}}{ \xi^{2k+4p} }\bigg)   = O (\frac{1}{\xi^{4p+i}} )\\
\frac{d^i}{d\xi^i} [ e^{\pm i 2(K-\xi)\ell }] &= & \frac{d^i}{d\xi^i} \big{[}
\sum_{n=0}^{\infty}\frac{(\pm i  2(K-\xi)\ell]^n}{n !}\big{]}=
\sum_{n=0}^{\infty}\frac{ (\pm i 2\ell)^n}{n !} \frac{d^i}{d\xi^i}(K-\xi)^n \\
\quad  \; & =& \left\{
\begin{array}{ll}
O(\frac{1}{\xi^{j+1}}),& j>0 \\
O(1)& j=0
\end{array}\right.
\end{array}\right.
$$
where we estimated 
$$\frac{d^j}{d\xi^j}(K-\xi)^n =\frac{d^j}{d\xi^j} \left(\frac{\eps^{2n}}{\xi^n}
\big[ \sum_{k=1}^{\infty}{ {\frac{1}{2}} \choose k} (-1)^k
(\frac{\eps^2}{\xi^2})^{k-1} \big]^n \right)$$

$$= O(\frac{1}{\xi^{n+j} }), \;\; n>0.$$

We have 
\begin{align*}
&\frac{d^n}{d\xi^n}\big[\psi (\lambda \xi) (K^2\xi^2)^{-p} e^{\pm i
2\ell(K-\xi)} \big]\\
=& \sum_{\substack{i+j+k=n\\ k>0}} \lambda^i \frac{1}{\xi^{4p+j}}
 \frac{1}{\xi^{k+1}} +\sum_{\substack{i+j+k=n\\ k=0}} 
 \lambda^i \frac{1}{\xi^{4p+j}} O(1)\\
\leq& C_{n,p,\ell}\lambda^{4p+1+n} +C_n \lambda^{4p+n} \leq
 \lambda^{4p+n}.   
\end{align*}

Integration by parts yields
$$ {\int^+}' =C_{n,\eps} \frac{ \lam^{4p+n-1}}{|x-y\pm 2\ell |^n}.
$$

It follows that 
\begin{equation}\label{eq:I+N}
I^+_N\leq C_{n,\eps} \sum_{p=0}^{N-1} \sum_{\ell=0}^{2p} \frac{ \lambda^{4p+n-1}}
{|x-y\pm 2\ell |^n}.
\end{equation}

Also, 
\begin{equation}\label{eq:R+N}
R^+_N\leq C \frac{ \lambda^{-1}}{|x-y |^n} \lam^{4N} \leq C \frac{
 \lambda^{n-1}}{|x-y |^n} \quad (4N\geq n)
\end{equation}
follows via integration by parts and the estimates
$$
\left\{
\begin{array}{ll}
\frac{d^i}{d\xi^i}[\psi (\lambda \xi)] \leq C {\lambda}^i \\
\frac{d^j}{d\xi^j}[(-1)^N (\frac{\eps^4 \sin^2 2K}{4K^2 {\xi}^2} )^N
\big/ (1+\frac{ \eps^4 \sin^2 2K }{ 4K^2 {\xi}^2 } ) ] =O(\frac{1}{ \xi^{4N}}). 
\end{array}\right.
$$

Combining (\ref{eq:I+N}) and (\ref{eq:R+N}) we obtain
$$
|I^+| \leq I^+_N + R^+_N \le  C_{N,\eps} \sum_{\ell=0}^{2N-2} \frac{\lam^{n-1}}{|x-y\pm 2\ell |^n}.$$

For $I^-$, we denote $\int^-= \int_{-1/\lam}^{-1/2\lam},$ then 
\begin{align*}
 I^- =&\int_{-1/\lam}^{-1/2\lam} \psi (\lam \xi) e^{i(x-y)\xi} d\xi
+\int^- \psi (\lam \xi)\overline{ {C'}_- } e^{i(x+y)\xi} d\xi\\
+& \int^- \psi (\lam \xi) C'_{-} e^{-i(x+y) \xi } d\xi
+\int^- \psi (\lam \xi) |C'_{-}|^2 e^{-i(x-y)\xi} d\xi\\
:=&I^-_1 +  I^-_2 + I^-_3+   I^-_4. 
\end{align*}






As  estimating $I^+$ we have 
$$ |I^-| \le C \sum_{\ell=0}^{2N}\frac{\lam^{n-1}}{|x-y\pm 2\ell |^n}. 
$$

Hence we obtain that if $x>1,\;y>1$,  
$$ 2\pi \vert K_j(x,y) \vert \le \vert I^+ \vert + \vert I^- \vert
\le C \sum_{\ell=0}^{2N} \frac{\lam^{n-1}}{|x-y\pm 2\ell |^n}.$$








\nd
$Case\; 2b. \;|x|\le 1, |y|\le 1 $.  Let notation be the same as in
Case 1a. 
by symmetry it is enough to deal with $ I^+$. 

From the expression of $B_+, B'_+$ we have 
\begin{align*}
I^+ &=\int^+ \psi(\lam \xi) (B_+ e^{iK x}+ B'_+ e^{-i K x})\overline
{B_+ e^{iK y}+ B'_+ e^{-i K y} } d\xi \\
\;& =\int^+ \psi(\lam \xi) |B_+|^2 e^{iK (x-y)} d\xi +\int^+ \psi(\lam
\xi) B_+ \overline{ {B'}_+} e^{iK (x+y)} d\xi \\
\;& + \int^+ \psi(\lam \xi)
B'\overline{B} e^{-iK (x+y)} d\xi +
\int^+ \psi(\lam \xi)|B'_+|^2 e^{-iK (x-y)} d\xi \\
\;&:=I_1^+ +I_2^+ + I_3^+ +I_4^+.
\end{align*}

We estimate  these terms separately. 
For instance, 
$$ I_2^+ =\frac{1}{4} \int^+ \psi(\lam \xi) e^{i (x+y)K} |C_+|^2
e^{-2 iK} (1-\xi^2/K^2) d\xi $$
Using the identity 
$$ |C_+ |^2 = \frac{ 4K^2 \xi^2}{ 4K^2 \xi^2 + \eps^4 \sin^2
2K}=\sum_{p=0}^{\infty}(-1)^p \left(\frac{ \eps^2 \sin 2K}{2K \xi} \right)^{2p}$$

$$=\sum_{p=0}^{N-1}  (-1)^p \left(\frac{ \eps^2 \sin 2K}{2K \xi} \right)^{2p}
+\widetilde{O} (\xi^{-4N}),$$
we can write 
$$ 4 I_2^+ = \sum_{p=0}^{N-1}  (-1)^p \int^+ \psi(\lam \xi) e^{i (x+y-2)K} 
  (\frac{ \eps^2 \sin 2K}{2K \xi} )^{2p} (1- \xi^2/K^2) d\xi $$

$$ + \int^+  \psi(\lam \xi) e^{i (x+y-2)K}\tilde{O} (\xi^{-4N}) (1-
\xi^2/{ K^2 } ) d\xi := I_{2,N}^+ +
R_{2,N}^+.$$
The integral in each term of the sum $I_{2,N}^+ $ is bounded by a
linear combination of the form 
$$ \int^+  \psi(\lam \xi) e^{i (x+y-2)K} e^{\pm i 2K \ell} (2K\xi)^{-2p}
(1-\xi^2/K^2) d\xi.$$
Integration by parts gives us
$$
|I_2^+| \le C_N \sum_{\ell=0}^{2N-1} \frac{\lam^{n-1}}{|x +  y \pm 2\ell |^n}. 
$$
The other terms $I_1^{+}, I_3^{+},  I_4^{+} $ also verify the above inequality
(possibly with $x+y$ replaced by $ x-y $).  And so does $ I^+ $
and $I^{-} $. 

We have 
$$ \vert K_j(x,y) \vert \le C \sum_{\ell=0}^{2N-1} \frac{\lam^{n-1}}{|x
\pm  y \pm 2\ell |^n}\; \;\;(|x|\le 1, |y|\le 1).$$ 
This completes the proof of Lemma 3.1(a). \hfill$\Box$



\vspace{.162in}
\nd
{\bf Proof of Lemma 3.3 (a).\quad} 
Note that $ \frac{\partial}{\partial x}e(x,\xi) $ exist for all $\xi
\neq \pm \eps , 0$ and are uniformly bounded in $x\in \R$ and 
$\xi$ in any bounded set.
$$  \frac{\partial}{\partial x} K_j(x,y)= \frac{1}{2\pi} \int_\R \vphi
(\lam^2 \xi^2) \frac{\partial}{\partial x}e(x,\xi) \bar{e}(y,\xi)
d\xi.$$ 

\nd
{\em Case} 1. $ x>1,\; y\in \R $.   
\begin{align*}
\frac{\partial}{\partial x} K_j(x,y)=&\frac{\partial}{\partial x}
\int_{ 2\eps < 1/ {2\lam} \le |\xi| 
\le 1/\lam } \psi(\lam \xi) ( C e^{ix
\xi } + C' e^{-ix \xi} ) \bar{e}(y,\xi) d\xi \\
=& \int i\xi  \psi(\lam \xi)(Ce^{ix
\xi } - C' e^{-ix \xi}) \bar{e}(y,\xi)d\xi \\
=& i \lam^{-1} \int \delta(\lam \xi) (Ce^{ix
\xi} - C' e^{-ix \xi}) \bar{e}(y,\xi)d\xi, 
\end{align*}
where $\delta(x)= x\psi(x) $ satisfies the same conditions as
$\psi(x)$:  $(i) \; \; \delta \in C^{\infty} \;\;\; (ii) \;\;\textrm{ supp}\;
\delta \subset \{\frac{1}{2} \le |\xi| \le 1 \} $ (except for $\psi$ being even, which is unimportant).   Thus we obtain, similar
to the case for $\psi(\lam H)(x,y)$ 
$$ 
\big\vert \frac{\partial}{\partial x} \int \big\vert \le C_N \sum_{\ell
=0}^{2N} \frac{\lam^{-2}}{ (1+\lam^{-1} |x\pm y\pm 2\ell |)^n}.
$$

\nd
\ca 2. $|x|\le 1,\; y\in \R $. 

Write $$ B(\xi)= \left\{
\begin{array}{ll}
B_+ & \xi >0\\
B_- & \xi <0.  
\end{array}\right.
$$

\begin{align*}
 \frac{\partial}{\partial x}\int =&  \frac{\partial}{\partial x}\int_{
 2\eps < 1/ {2\lam} \le |\xi| \le 1/\lam } \psi(\lam \xi)(B e^{iK x}
 + B' e^{-i K x}) \bar{e}(y,\xi)d\xi \\
=&\int_R iK
 \psi(\lam \xi)(B e^{iK x}
 - B' e^{-i K x}) \bar{e}(y,\xi)d\xi \\
=& i \lam^{-1} \int (\lam \xi) \psi(\lam \xi) (B e^{iK x}
 - B' e^{-i K x}) K/\xi \;\overline{e}(y,\xi)d\xi \\
=& i \lam^{-1} \int  \delta(\lam \xi) (B e^{iK x}
 - B' e^{-i K x}) K/\xi \; \overline{e}(y,\xi) d\xi,
\end{align*}
where $\delta(x)= x \psi(x) $.
Thus we obtain, similar to the case for $\psi(\lam H)(x,y)$ 
$$ 
\big\vert \frac{\partial}{\partial x} \int \big\vert \le C_N \sum_{\ell
=0}^{2N} \frac{\lam^{-2}}{ (1+\lam^{-1} |x\pm y\pm 2\ell |)^n}.
$$

\nd 
\ca 3. $ x<-1,\; y\in \R $. 
The relation $\vphi(\lam H)(x,y)= \vphi(\lam H)(-x,-y)$ implies 
\begin{align*}
&\frac{\partial}{\partial x} [\vphi(\lam H)(x,y)]_{x<-1}=
\frac{\partial}{\partial x} [\vphi(\lam H)(-x,-y)]_{x<-1} \\
=& -[ \frac{\partial}{\partial x} \phi(\lam H)] (-x,-y)]_{x<-1}.
\end{align*}
Therefore, if $x<-1$
\begin{align*}
&\big\vert \frac{\partial}{\partial x} \vphi(\lam H)(x,y) \big\vert =\big\vert
\frac{\partial}{\partial x} \big[\vphi(\lam H)(x,y)\big] (-x, -y) \big\vert\\
\le &C
\sum_{\ell =0}^{2N} \frac{\lam^{-2} }{ (1+\lam^{-1} |- x \mp y\pm 2\ell
|)^n}.
\end{align*}
This concludes the proof of Lemma 3.3(a).\hB

\vspace{.2in}
4.2. {\em Low energy estimates}$\quad j\le J$

\vspace{.162in}
4.2.1. {\bf Proof of Lemma 3.1(b). \quad} 
We study the decay of the kernel of $\vphi_j(H)$ as $ j \rightarrow -\infty$. 
As in the high energy case, we only need to check four cases 
 1a, 1b, 1c and 2b. Outline will be given for 1a, 2b only.

\nd
\ca 1a. $x>1,\; y>1$  ($1/\lam \le4\eps$)  
$$
2\pi K_j (x,y)= \int_\R \psi(\lam \xi) e(x,\xi)\bar{e}(y,\xi)d\xi
= \int^+ +\int^- .$$
We obtain by integration by parts
\begin{align*}
\big{\vert} \int^+ \big\vert 
\le& C \frac{\lam^{n-1}}{|x-y|^n},
\end{align*}
where we used 
\[
\left\{
\begin{array}{ll}
\frac{d^n}{d\xi^n }[\psi(\lam \xi)] = \lam^n \psi^{(n)}(\lam \xi) \le
C \lam^n,  \quad \lam =2^{-j/2}\rightarrow +\infty \; (j\rightarrow
-\infty) \\
\frac{d^i}{d\xi^i }(|C_+|^2) =\left\{
\begin{array}{ll}
O(\xi^2)& i=0\\
O(\xi) & i=1\\
O(1) & i>1
\end{array}\right.
\end{array}\right. 
\]


We obtain also
\begin{align*}
 \vert \int^- \vert \le& C_n \frac{\lam^{n-1} }{|x-y|^n},
\end{align*}
using 
\begin{align*}
\frac{d^i}{d\xi^i}\big( |C'_-|^2 \big) =&
\frac{d^i}{d\xi^i} \big[\frac{ \eps^4 }{ (2\rho/ {\sh 2\rho} )^2 \xi^2 +
\eps^4 } \big] = O(1), \quad \xi \rightarrow 0^-\\
\textrm{and}\;  \frac{d^i}{d\xi^i} {C_-}' =&\;O(1),\; \textrm{since}\; C'_- \;\textrm{is}\;
C^{\infty}_{[-4\eps, 0)}.   
\end{align*}







\nd
\ca 2b. $ |x|\le 1,\; |y|\le 1.$  Let $I^+(x, y):=\int^+$, 
$I^-(x, y):=\int^-$. 
\begin{align*}
 &I^+(x, y) = \int^+ \psi( \lam \xi) (B_+ e^{\rho x} +B'_+ e^{-\rho x})
\overline{ B_+ e^{\rho y} +B'_+ e^{-\rho y} } d\xi \\
= &\int^+ |C_+|^2 ( \cosh \rho (1-x) -i \xi/\rho \sinh \rho (1-x) ) 
( \cosh \rho (1-y) +i \xi/\rho \sinh \rho (1-y) ) d\xi\\
 \le &C \lam^{-1} \le 3^n C \frac{\lam^{-1} }{1+ \lam^{-1}|x \pm
y|)^n}, 
\end{align*}
where we note that $\cosh \rho (1-x) -i \xi/\rho \sinh \rho (1-x)$ and $\cosh
\rho (1-y) +i \xi/\rho \sinh \rho (1-y)$ are bounded by a constant
uniformly in $|x| \le 1$ and $|y| \le 1$.

The term  $I^-(x, y)$ satisfies the same inequality by the relation $  
I^-(x, y)= I^+(-x, -y)$.
\hfill $\Box$






\vspace{0.2in}
\nd
{\bf Proof of Lemma 3.3 (b).\quad}  
The same argument in proving Lemma 3.1(b) is valid for the proof of 
Lemma 3.3(b). The interested reader can fill in the details.

\section{Local energy estimates}

Let $\Phi \in C^{\infty}$ have support contained in $\{\xi: |\xi| \le 1
\}$. Then the spectrum of $ \Phi(H)$ includes the low energy, a
neighborhood of $0$. We use the term ``local energy'' to distinguish
from the low energy case, where the support of $\vphi_j (j\le J)$ keeps away
from 0. Since $0\in \textrm{supp} \;\Phi$ and $e(x,\xi)$ is not
continuous at the origin $\xi=0$ (!), we need to treat the corresponding
kernel more carefully. The proof is more delicate and requires a
``matching'' method.   

\vspace{.15in}
\nd
{\bf Proof of Lemma 3.2.\quad} As in $\S$4, the estimates rely on 
four cases 1a, 1b, 1c, 2b. 
We use $\hat{f}$ and $\check{f}$ to denote the 
ordinary Fourier transform and its inverse, resp. 

\nd
\emph{Case} 1a. $x>1, y>1$.  Let $K(x,y)= \Phi(H)(x,y),  \;\Psi(x)=\Phi(x^2)$.  
\begin{align*}
2\pi K(x,y)=& \int_0^1 \Psi(\xi) C_+ e^{ix \xi} \overline{ C_+ e^{i y \xi}
} d\xi + \int^0_{-1} \Psi(\xi) (e^{ix \xi} +C'_- e^{-ix \xi}) \overline{
e^{iy \xi} + C'_- e^{-iy \xi } } d\xi \\
 =& I^+ + I^-.
\end{align*}
Write 
\begin{align*}
I^- = & \int^0_{-1} \Psi(\xi) e^{i(x-y) \xi }d\xi +
\int^0_{-1}\Psi(\xi) C'_- e^{-i(x+y) \xi} d\xi \\
+& \int^0_{-1} \Psi(\xi)
\overline{C'_-} e^{i(x+y)\xi}d\xi +   \int^0_{-1} \Psi(\xi) |C'_-|^2
e^{-i(x-y )\xi}  d\xi \\
 :=&  I_1^- + I_2^- +  I_3^- +  I_4^-.
\end{align*}
The relations $ C'_- (-\xi) =   A'_+ (\xi)= \overline{ C'_-} (\xi)$ and 
$ |C_+|^2 +   |A'_+|^2 =  |C_+|^2 +   |C'_-|^2=1 $ imply that 
\begin{align*}
&I^+ + I_1^- + I_4^- = \int^1_0  \Psi(\xi) |C_+|^2 e^{i(x-y) \xi }
d\xi \\
+& \int^1_0  \Psi(\xi) e^{i(x-y) \xi }
d\xi + \int^1_0  \Psi(\xi) |C'_-|^2 e^{i(x-y) \xi }
d\xi  \\
 =& \int^1_{-1} \Psi (\xi)  e^{ i (x-y) \xi }
d\xi  = \sqrt{2\pi} {\Psi }^\vee (x-y). 
\end{align*}  

Also, the relation $ C'_- (-\xi) = \overline{C'_-} (\xi)$ gives 
\begin{align*}
I_2^- + I_3^- =& \int_{-1}^0 \Psi(\xi) C'_- e^{-i(x+y) \xi} d\xi 
+  \int_{0}^1 \Psi(\xi) C'_- (\xi) e^{-i(x+y) \xi} d\xi \\
=& \sqrt{2\pi} (\Psi(\xi) C'_-(\xi))^{\wedge} (x+y).
\end{align*}

Since $ \Psi\in C_0^{\infty}$ and $ C'_- \in C^{\infty}$, we have for
$ x>1, y>1$ 
\begin{align*} 
 &2\pi \vert K(x,y) \vert \le \vert  I^+ + I_1^- + I_4^- \vert +  \vert
 I_2^- + I_3^- \vert  \\
 \le& \frac{C_n}{ (1+ |x-y|)^n} + \frac{C_n}{ (1+ |x+y|)^n } \le \frac{C_n}{ (1+ |x-y|)^n}
\end{align*}
by the rapid decay for the Fourier transform of
$C_0^{\infty}$ functions, where 
$$
C'_- (\xi) =i \frac{\eps^2\sinh 2\rho}{2\rho \xi} A_-(\xi) = i
\frac{\eps^2 e^{2i\xi} \sinh 2\rho / 2\rho} {\xi \cosh 2\rho - i
(\rho^2-\xi^2)\sinh 2\rho / 2\rho }  \in C^{\infty} (\R).$$

\nd
\ca 1b. $x>1, |y|\le 1$

Using $ e_+(y, \xi)= C_+ e^{i\xi} [ \cosh \rho(1-y)-i  \xi/\rho \sinh
\rho(1-y)]$ and $ A_- = \frac{  2\rho \xi e^{2 i\xi} }{  2\rho \xi
\cosh 2\rho -i (\rho^2-\xi^2)\sinh 2\rho} $,

\begin{align*}
&2\pi  K(x,y)= \int_0^1 \Psi (\xi)  e^{i(x-1) \xi} |C_+|^2 \underbrace{\Re
+i \Im}  [\cosh \rho (1-y) + i \xi/\rho \sinh \rho (1-y)] d\xi \\
 +& \int_{-1}^0  \Psi(\xi)  e^{i(x-1) \xi} \underbrace{\Re +i
\Im} \frac{ 2\rho \xi[\cosh \rho (1+y)-i  \xi/\rho \sinh
\rho (1+y)] }{ 2\rho \xi \cosh 2\rho +i (\rho^2-\xi^2)\sinh 2\rho}
+ \\
  &\int_{-1}^0  \Psi(\xi)  e^{-i(x-1) \xi} \underbrace{\Re +i
\Im} \frac{ i\eps^2 \sinh 2\rho \cdot 2\rho \xi}{ 4\rho^2 \xi^2
+\eps^4 \sinh^2 2\rho}   [\cosh \rho (1+y)-i  \xi/\rho \sinh
\rho (1+y)] d\xi\\
:=& ``{ \cal R}e\text{ ''} + i `` {\cal I}m \text{''},
\end{align*}
where we break each of the above three integrals into two parts; then
let 
$``\mathcal{ R}e \text{''}$ be the sum of the three integrals involving real parts only,
and let $``\mathcal{ I}m \text{''}$ be the sum  of the three integrals involving
imaginary parts only. 

We have 
\begin{align*}
`` \mathcal{R}e\text{''}= \int_0^1&  \Psi(\xi) |C_+|^2 e^{i(x-1) \xi}\cosh \rho
(1-y)d\xi+ \\
\int_{-1}^0&  \Psi(\xi) e^{i(x-1) \xi} \Re\big[\frac{ 2\rho \xi (\cosh \rho
(1+y)-i  \xi/\rho \sinh
\rho (1+y) ) }{ 2\rho \xi \cosh 2\rho +i (\rho^2-\xi^2)\sinh 2\rho} \big]
d\xi \\
 + \int_{-1}^0&  \Psi(\xi) e^{-i(x-1) \xi}\Re\big[ C'_- \overline{A_-} (\cosh
\rho (1+y)-i  \xi/\rho \sinh\rho (1+y) )\big] d\xi\\
=\int_0^1& \Psi(\xi) e^{i(x-1) \xi} \frac{4\rho^2 \xi^2 \cosh \rho
(1-y) + 2\eps^2 \xi^2   \sinh 2\rho \sinh \rho (1+y)}{ 4\rho^2 \xi^2
+\eps^4 \sinh^2 2\rho}+\\
 \int_{-1}^0  \Psi(\xi)& e^{i(x-1) \xi} \frac{4\rho^2 \xi^2 \cosh 2\rho
\cosh \rho (1+y) - 2 \xi^2(\rho^2 -\xi^2)   \sinh 2\rho \sinh \rho (1+y)}{ 4\rho^2 \xi^2 +\eps^4 \sinh^2 2\rho} d\xi. 
\end{align*}
Noting that $\rho^2 -\xi^2 = 2\rho^2 -\eps^2$, and \\
$ \ch 2\rho \cosh 
\rho (1+y) -  \sh 2\rho \sinh
\rho (1+y) = 
\cosh \rho (1-y)$
we obtain
\begin{equation}\label{eq:Re}
`` \mathcal{R}e\text{''}= \sqrt{2\pi} { \big[  \Psi(\xi) e^{-i \xi} \frac{4\rho^2 \xi^2 
\cosh \rho (1-y) + 2 \eps^2 \xi^2 \ch \rho(1-y)+ 2 \eps^2 \xi^2   \sinh 2\rho \sinh \rho(1+y)}{ 4\rho^2 \xi^2
+\eps^4 \sinh^2 2\rho} \big] }^\vee (x). 
\end{equation} 

For 
\begin{align*}
 ``\mathcal{I}m \text{''} &= \int_0^1   \Psi(\xi) |C_+|^2 e^{i(x-1) \xi} \xi/ \rho \sh\rho (1-y)d\xi\\
+\int_{-1}^0  \Psi(\xi)& e^{i(x-1) \xi} \Im\big\{ \frac{ 2\rho \xi [\cosh
\rho (1+y)-i  \xi/\rho \sinh
\rho (1+y) ] }{ 2\rho \xi \cosh 2\rho +i (\rho^2-\xi^2)\sinh 2\rho} \big\}
d\xi\\
+ \int_{-1}^0  \Psi(\xi)& e^{-i(x-1) \xi} \frac{\eps^2 \sh 2\rho 2
\rho \xi}{ 4\rho^2 \xi^2
+\eps^4 \sinh^2 2\rho}   \cosh \rho (1+y) d\xi\\
= \int_0^1 \Psi(\xi) &e^{i(x-1) \xi} \frac{ 2
\rho \xi}{ 4\rho^2 \xi^2
+\eps^4 \sinh^2 2\rho}  [ 2\xi^2 \sh \rho(1-y) - \eps^2 \sh 2\rho
\ch\rho (1+y)] d\xi +\\
 \int_{-1}^0 \Psi(\xi) e^{i(x-1) \xi}& \frac{ 2
\rho \xi}{ 4\rho^2 \xi^2
+\eps^4 \sinh^2 2\rho}  [- 2\xi^2 \ch 2\rho \sh \rho(1+y) - (\rho^2- \xi^2) \sh 2\rho\ch\rho (1+y)] d\xi
\end{align*}
Noting that $\rho^2 -\xi^2 = \eps^2  -2\xi^2$, and\\
 $ \sh 2\rho \cosh 
\rho(1+y) -  \ch 2\rho \sinh
\rho(1+y) = 
\sh \rho(1-y)$
we obtain
\begin{equation}\label{eq:Im}
``\cal{I}m\text{''}= \sqrt{2\pi} \big[  \Psi(\xi) e^{-i \xi} \frac{2\rho
\xi}{4\rho^2 \xi^2
+\eps^4 \sinh^2 2\rho}  ( 2\xi^2 \sh \rho (1-y) - \eps^2 \sh 2\rho
\ch\rho (1+y) ) \big]^\vee (x).
\end{equation}
Since the functions in the square brackets of (\ref{eq:Re}) and 
(\ref{eq:Im})
 are $C^{\infty} $, it follows that for $x>1, |y|\le 1$,  
$$
\vert K(x,y)\vert \le \frac{C_n}{(1+|x|)^n} \le \frac{C'_n}{(1+|x-y|)^n },$$
where $C_n$ can be taken to be independent of $|y|\le 1$.

\nd
\ca 1c. $\; x>1, y< -1$. The proof is similar to that of Case 1a and hence
omitted.




\nd
\ca 2b. $ |x|\le 1, |y|\le 1$. Since $ |e(x,\xi)| \le C_{\eps} $, for
all $x, \xi$,  the result is trivial:
$$ \vert K(x,y) \vert \le C'_{\eps} \sim \frac{C_n}{ (1+|x-y|)^n} $$
whenever $  |x|\le 1, |y|\le 1 $. This concludes the proof of Lemma
3.2. $\hfill \Box$

\vspace{.2in}
\nd
{\bf Proof of Lemma 3.4.\quad} 
 With the convention $ \int^+ =\int_0^1, \; \int^- =\int_{-1}^0$, 
\begin{align*}
&2\pi \frac{\partial}{\partial x}  K(x,y) = \frac{\partial}{\partial x}
\int_{-1}^1 \Psi(\xi) e(x,\xi) \bar{e}(y,\xi)d\xi\\
:=&\frac{\partial}{\partial x}  \int^+ + \frac{\partial}{\partial x}
\int^- .
\end{align*}
The function $ \xi \mapsto \frac{\partial}{\partial x}  e(x,\xi) $
is discontinuous at $ \xi =0$. As suggested by  the treatment of $
K(x,y) $ we want to properly ``match'' different parts of the above
integrals so $ \frac{\partial}{\partial x} K(x,y) $ can be written as
a linear combination of the Fourier transform of $
C_0^{\infty} $ functions.

We only need to check five cases $1a, 1b, 1c, 2a, 2b$. Estimates for
the other cases follow readily from the relation $
\frac{\partial}{\partial x} K(x,y)=  \frac{\partial}{\partial x}[
K(-x,-y )]
= - (\frac{\partial}{\partial x}
K) (-x,-y ).$  We outline the proofs for 1a, 1b and 2b only, since   
1c and 2a can be dealt with similarly.

\vspace{.152in}
\nd
\ca 1a. $x>1, y>1$. Let $\Delta(\xi)= i \xi \Psi(\xi) \in C_0^{\infty} $.
\begin{align*}
 \frac{\partial}{\partial x}  \int^+ =& \int^+ \Psi(\xi)i \xi |C_+|^2 e^{i(x-y)\xi}
d\xi =  \int^+ \Delta(\xi) |C_+|^2 e^{i(x-y)\xi} d\xi .\\
 \frac{\partial}{\partial x}  \int^- =& \int^- i\xi \Psi(\xi) (e^{ix
\xi}- C'_- e^{-ix \xi}) \overline{e^{iy
\xi}+ C'_- e^{-iy \xi}  }d\xi\\
 =&  \int^- \Delta(\xi)  e^{i(x-y)\xi} d\xi - \int^- \Delta(\xi)
|C'_-|^2 e^{-i(x-y)\xi} d\xi\\
+&  \int^- \Delta(\xi) \overline{C'_-} e^{i(x+y)\xi} d\xi - \int^- \Delta(\xi) C'_- e^{-i(x+y)\xi} d\xi\\
= &\int^-  \Delta(\xi)  e^{i(x-y)\xi} d\xi + \int^+ \Delta(\xi)
|C'_-|^2 e^{i(x-y)\xi} d\xi\\
 -& \int^+ \Delta(\xi)
C'_-  e^{-i(x+y)\xi} d\xi - \int^- \Delta(\xi)
C'_-  e^{-i(x+y)\xi} d\xi , 
\end{align*}
where we note $ \Delta(\xi) $ is odd and the relation $ C'_- (\xi) = \overline{ C'_-} (-\xi).$\\
We have, by the relation $ |C_+|^2 + 
|C'_-|^2 =1$,
\begin{align*}
\frac{\partial}{\partial x}  \int^+ &+\frac{\partial}{\partial x}
\int^-  = \int   \Delta (\xi) e^{i(x-y)
\xi}d\xi - \int  \Delta (\xi) C'_- e^{-i(x+y) \xi} d\xi \\
=& \sqrt{2\pi} { [\Delta(\xi)] }^\vee (x-y) - \sqrt{2\pi} {
[\Delta(\xi) C'_- ] }^\wedge (x+y).
\end{align*}
Since $\Delta \in C^{\infty}_0, C'_- (\xi) \in C^{\infty} $, 
 the inequality in Lemma 3.4  holds for $ x>1, y>1$. 

\nd
\ca 1b.  $ x>1, |y| \le 1$. Let notation be as in Case 1a.
\begin{align*}
 & \frac{\partial}{\partial x}  \int^+ = \int^+  \Delta(\xi)
|C_+|^2 e^{i(x-1)\xi} (\Re +i \Im) [ \ch \rho (1-y) + i \xi/ \rho \sh
\rho (1-y) ]   d\xi\\
&  \frac{\partial}{\partial x}  \int^- = \int^-  \Delta(\xi)
 e^{i(x-1)\xi} (\Re +i \Im) \big[ \frac{ 2\rho \xi (\ch \rho (1+y) - i \xi/ \rho\; \sh\rho (1+y)}{ 2\rho \xi \ch \rho + i(\rho^2 - { \xi}^2 )  \sh
2\rho }   \big]   d\xi\\
- \int^- &\Delta(\xi)
 e^{-i(x-1)\xi} (\Re +i \Im) \big[\frac{i \eps^2 \sh 2\rho \cdot 2\rho \xi}{
4 \rho^2 \xi^2 + \eps^4 \sh^2 2\rho} (\ch \rho (1+y) - i \xi/ \rho \sh
\rho (1+y) ) \big]   d\xi .
 \end{align*}

As in the case for $K(x,y)$, we split each integral into two parts
and let ``$\cal{R}e$'' and ``$\cal{I}m$'' denote the sum of integrals involving only
these ``reals''  and ``imaginaries'' respectively.  As a result ,
$$
2\pi \frac{\partial}{\partial x} K(x,y) = ``\mathcal{R}e\text{''} + i ``\mathcal{I}m\text{''},
$$
where we find, by noting that $ \Delta$ is odd,
$``\mathcal{R}e\text{''}$  and $``\mathcal{I}
m\text{''}$
have the same expressions as in (\ref{eq:Re}) and (\ref{eq:Im}) resp., 
except that $\Psi$ should be replaced by $\Delta$. 
Case 1b is so verified. 

\vspace{.12in}







Finally, the decay estimate for Case 2b ($ |x|, |y| \le 1) $ follows trivially from the
fact that $e(y,\xi) \in L^{\infty} (\R \times \R )$ and $ \frac{\partial}
{\partial x} e(x,\xi) \in L^{\infty}_{loc} (\R \times
\R )$ 

$$\vert  2\pi K(x,y)\vert = \Big\vert \int_{|\xi| \le 1} \Psi(\xi)
\frac{\partial}{\partial x} e(x,\xi) \bar{e}(y,\xi) d\xi \Big\vert $$ 

$$ \le C \le  \frac{C'_n }{ (1+|x-y|)^n},$$
where 
$$ \frac{\partial}{\partial x} e_+ (x,\xi) = C_+ e^{i\xi} [ - \rho  
 \sh \rho (1-x) + i\xi   \ch  \rho (1-x)],  $$ 

$$  \frac{\partial}{\partial x} e_- (x,\xi) = A_- e^{-i\xi} [  \rho  
 \sh \rho (1+x) + i\xi   \ch  \rho (1+x) ] . $$
This completes the proof of Lemma 3.4.
\hfill $\Box$

\section{Spectral multipliers for $L^p$ and $B_p^{\alpha , q}(H)$ }
The operator $m(H)$ can be defined using functional calculus: $m(H)=
\int_0^{\infty} m(\lam) d E(\lam) $, where ${H} = \int_0^{\infty}
\lam d E(\lam) $  is the spectral resolution of $H$. 
From \cite{Zhen, Zhen2}, we know $ m(H)$ also has the expression $ m(H)=
{\mathcal{F}}^{-1} m(\xi^2)\cal{F}$ if $m$ is bounded.

We shall prove that under the same differentiability condition on 
$m\in L^\infty$ as in
the Fourier case, $m(H) $ has a bounded extension on $L^p$ from $L^p
\cap L^2$, by showing that the 
kernel $m(H)(x,y)$
satisfies a H\"{o}rmander type condition:  (compare the Fourier case
[Stein 93]):
\begin{equation}
\int_{ z>2|y-\bar{y} | } |m(H)(x,y)- m(H)(x,\overline{y })| dx \le A,
\end{equation}
where $ z=\min_\pm( |x  \pm \overline{y} |)$. 

An immediate question arises: what is the kernel expression for
$m(H)$?  Since $m$ may not necessarily have compact support, the
answer is not so immediate. Let $\{\de_j\}_{-\infty}^\infty$ be a smooth dyadic resolution
of unit and $m_j(x)=m\de_j(x)$. Then for $ f\in
L^2, \; m(H)f = \sum_{-\infty}^{\infty} m_j(H)f $ in $L^2$. This
suggests that $ m(H)(x,y) $ may have the (pointwise) expression $
\sum_{-\infty}^{\infty}  K_j(x,y)$, where
$K_j$ denotes the kernel of $m_j(H)$.
Our next lemma shows that this is true in an appropriate sense.

\begin{lemma} \label{lem:m-ker}
Let $m$ be bounded and $ |m'(\xi) | \le C |\xi|^{-1} $ for \mbox{$ \xi \in
\R\setminus \{ 0 \} $.}
Then for $ f\in L^2_0= \{ f\in L^2: \textrm{supp}\; f\;\textrm{is} \;
compact\}$, $m(H)f $ has
the expression 
$$ m(H)f(x)= \int  m(H)(x,y) f(y)dy $$
for $a.e.  \; x \notin \pm \textrm{supp} \;f$, where $  m(H)(x,y) = \sum_{-\infty}^{\infty} m_j(H)(x,y).$ 
\end{lemma}

\nd
\pf Since $\sum_{-\infty}^{\infty} m_j(H)f$ converges to 
$ m(H)f$ in $L^2$,    
 it suffices to show the series $ \sum m_j(H)f(x) $ converges pointwise
for each $ x \notin \pm
\textrm{supp} f$. 

 Let $ 0<t=$  the distance from $x$ to the set
$ (\textrm{supp} \;f )\cup (-\textrm{supp} \;f)$. Then $\textrm{supp}\; f \subset \{ y: \min (|y+x|, |y-x |) \ge t\}$.   
By Lemma \ref{lem:ker-size} 
we have for  $ x \notin \pm \textrm{supp} \;f $, picking $J\in \Z$ 
\begin{align*}
&\sum_{-\infty}^J \Big\vert \int m_j(H)(x,y)f(y) dy \Big\vert \leq \Vert f
\Vert_2  \sum_{-\infty}^J  \Vert  m_j(H)(x,\cdot) \Vert_2 \\
\le& C   \Vert f
\Vert_2  \sum_{-\infty}^J 2^{j/4} \le C_J \Vert f \Vert_2,
\end{align*}    
and writing $\min |y\pm x|=\min(|y+ x|,|y-x|)$
\begin{align*}
 &\sum_{J+1}^{\infty} \Big\vert \int m_j(H)(x,y)f(y) dy \Big\vert =
\sum_{J+1}^{\infty} \Big\vert \int_{ \min |y\pm x| >t}  m_j(H)(x,y)f(y)
dy \Big\vert \\
\le &\Vert f \Vert_2  \sum_{J+1}^{\infty} ( \int_{\min |y\pm x| >t}
\vert  m_j(H)(x, y) \vert^ 2 dy )^{1/2} \\
=  & C \Vert f \Vert_2 t^{-1 } \sum_{J+1}^{\infty} 2^{-j/4}  
\le C_J  \Vert f \Vert_2 t^{-1 },
\end{align*}
where we used the inequality $( \int_{\min |y\pm x| >t}
\vert  m_j(H)(x, y) \vert^ 2 dy )^{1/2}\le Ct^{-1}2^{-j/4}  $,
by (\ref{eq:zker-size}).
This shows that $\sum m_j(H)f(x) $ converges for all $ x \notin \pm \textrm{supp} f $. 

\quad \hfill $ \Box$

\begin{lemma} \label{lem:ker-size} Let $ z= \min |x\pm y|$ and 
$\lam = 2^{-j/2}$. Then there exists a constant $C$ independent of $y$
so that
\begin{align} 
\Vert & K_j (\cdot , y)\Vert_2 \le C \lam^{-1/2}, 
\label{eq:m-K_j} \\
\Vert & z K_j (\cdot , y)\Vert_2 \le C \lam^{1/2}, 
  \label{eq:zker-size}\\
\int_{|z|>t} &\vert K_j (x , y)\vert dx \le C   
t^{-\frac{1}{2}}\lam^{1/2}. 
\end{align}
\end{lemma} 

\begin{lemma} \label{lem:der-ker-size} Let $z$, $\lam$ be as above. Then there
exists a constant $C$, independent of $y$ so that
\begin{align}
\Vert &\frac{\partial}{\partial y}  K_j (\cdot , y)\Vert_2 \le C \lam^{-3/2}, \; \label{eq:m-der-K}\\
\Vert & z \frac{\partial}{\partial y} K_j (\cdot , y)\Vert_2 \le C \lam^{-1/2}, \; 
 \label{eq:m-zder-K}\\
\int_{|z|>t}& \Big\vert \frac{\partial}{\partial y} K_j (x , y)\Big\vert dx \le
C   t^{-1/2}\lam^{-1/2}. \; 
 \label{eq:int-m-der-K}
\end{align}
\end{lemma}

We are ready to verify the H\"{o}rmander condition for $ m(H)$.

\begin{lemma} Let $z = \min |x\pm \overline{y} |, \;t=|y - \overline{y}|$
and $\lam =2^{-j/2}$.
Then
\begin{eqnarray}\label{eq:pre-Horm}
\int_{|z|>2t} \vert K_j (x , y)- K_j (x , \overline{y}) \vert  \ud x  \\
\nonumber
\le C \left\{
\begin{array}{lcl}
  t^{1/2}\lam^{-1/2} & if & t\lam^{-1} \le 1 \\
 t^{-1/2}\lam^{1/2} & if& t\lam^{-1} > 1 .
\end{array}\right.
\end{eqnarray}
Moreover,
\begin{equation}\label{eq:Horm} 
\int_{|z|>2t} \vert K (x , y)- K (x , \overline{y}) \vert  \ud x \le A,
\end{equation}
where $ K(x,y) $ agrees with a ``function'' in the sense of 
Lemma \ref{lem:m-ker}.
\end{lemma}

\nd \rk
Compare \cite{Ste3,E2} for Lemma 6.2--6.4.

\nd
\pf  Let $ y\in \overline{y} +I,  I =[-t,t]$.  If $ t\lam^{-1} \le 1$, by Lemma 6.3
\begin{align*}
 &\int_{ \{|x-\bar{y}| >2t \} \cap \{|x+\bar{y}| >2t \}  }
| m_j(H) (x , y)- m_j(H) (x , \bar{y}) \vert  dx \\
=& 
\int_{z>2t} \Big\vert
\int_{\bar{y}}^y  \frac{\partial}{\partial \xi}  m_j(H) (x , \xi) d\xi
\Big\vert dx \\
\le &\int_{\bar{y}}^y d\xi \int_{ z>2t}  \Big\vert \frac{\partial}{\partial
\xi}  m_j(H) (x , \xi)\Big\vert  dx \\
\le &t \sup_{ \xi \in \bar{y} +I} \int_{  \{|x-\bar{y}| >2t \} \cap \{|x+\bar{y}| >2t \} } 
 \Big\vert \frac{\partial}{\partial \xi}  m_j(H) (x , \xi) \Big\vert dx\\
\le &C_{\eps}  t^{1/2}\lam^{-1/2}, \; \textrm{all} \;\bar{y},
\end{align*}  
If $t\lam^{-1} >1$,  
\begin{align*}
& \int_{ |z| >2t } \vert  K_j (x , y)- K_j (x , \bar{y}) \vert dx \\
\le& \int_{ |z| >2t } \vert  K_j (x , y) \vert dx + \int_{ |z| >2t }
\vert K_j (x , \bar{y}) \vert dx\\
\le& \int_{ \min | x \pm y |> t } \vert  K_j (x , y) \vert dx + \int_{\min
|x\pm \bar{y} |>2t}  \vert  K_j (x , \bar{y}) \vert dx \\
\le& C t^{-1/2}\lam^{1/2}.  
\end{align*}
This proves (\ref{eq:pre-Horm}).

Now we show (\ref{eq:Horm}) using (\ref{eq:pre-Horm}). Let $ h= |y-\bar{y}|$. By Lemma 6.1,
\begin{align*}
 &\int_{z>2|y-\bar{y}| } | m(H)(x, y) - m(H)(x, \bar{y})| dx\\
\le& \sum_{-\infty}^{\infty} \int_{ z>2|y-\bar{y}| } | m_j (H)(x, y)
- m_j (H)(x, \bar{y})| dx :=  \sum_{ |h| \lam^{-1} \le 1} + \sum_{|h| \lam^{-1} > 1 }.  
\end{align*}
By Lemma 6.4 the first sum is bounded by 
$$ C h^{1/2} \sum_{ h 2^{j/2} \le 1}  2^{j/4} \le C $$
and the second sum by 
\begin{align*}
 &\sum_{ h 2^{j/2} > 1} \Big( \int_{  z>2 h } | m_j (H)(x, y)|dx
+ \int_{  z>2 h} |m_j (H)(x, \bar{y})| dx \Big)  \\
\le &  \sum_{ h  2^{j/2} > 1} h^{-1/2}
2^{-j/4} \le C,
\end{align*}
where we note that $ \int_{  z>2 h } | m_j (H)(x, y)|dx \le \int_{
\min | x\pm y |> h } | m_j (H)(x, y)|dx$, \\
because  $|z|>2 h$ implies that $ |x\pm y | \ge
|x\pm \bar{y}| - |y-\bar{y}|> 2h-h=h.$  Hence (\ref{eq:Horm}) holds. 
\quad\hfill $\Box$

\begin{theorem}\label{th:L^p-mult} Suppose $m\in L^{\infty}: \R \rightarrow \mathbf{C}$ satisfies 
$|m'(x)|\leq C |x|^{-1}$. Then $m(H)$ is bounded on $L^p, 1<p<\infty$
and of weak type $(1,1)$.
\end{theorem}

As a consequence of Theorem 6.5, we shall show that, $m(H)$ initially
defined for $f\in L^2 $, has  a bounded linear extension to the Banach
spaces $ B_p^{\alpha , q} (H), 1<p<\infty$.

\begin{theorem} \label{th:Besov-mult} Suppose $m\in L^{\infty}$ be as above. Then $m(H)$ extends
to a bounded linear operator on $B_p^{\alpha,q}(H) $ for $ 1<p<\infty,
0<q\leq \infty, \alpha \in \R$. 
\end{theorem}

\vspace{.17in}
\nd
{\bf Proof of Theorem \ref{th:L^p-mult}.}  Applying Calder\'on-Zygmund 
decomposition and using Lemma 6.4,
we can get the $weak\; (1,1)$ result for $m(H)$. Then the $L^p $ result,
$1<p<\infty $, follows via Marcinkiewicz interpolation and duality.
For completeness, we give the $weak\; (1,1)$ estimation. It is enough
to assume $f\in L^1\cap L^2$ by density.

Given $f\in L^1, s>0$, according to Calder\'on-Zygmund lemma 
there is a decomposition $f= g+ b $ with $ b=\sum b_k$ and a countable
collection of disjoint open intervals $ I_k$ such that 

i) $|g(x) | \le C s $ a.e.    
  
ii) Each $b_k $ is supported in $I_k, \int b_k dx =0$ and
$$ s \le \frac{1}{ |I_k| } \int_{I_k} |b_k| dx \le 2s $$

iii)  Let $ D_s ={\cup}_k I_k = {\cup}_k (\bar{y}_k- t_k,  \bar{y}_k+ t_k)$, 
where $ 2t_k = |I_k|>0, \bar{y}_k $ is the center of $ I_k$.
Then 
$$ |D_s| \le C s^{-1} \Vert f\Vert_1.$$

iv) $ g\in L^1 \cap L^2$, $g(x) = f(x) \; \textrm{if}\; x\notin D_s$, and 
\begin{equation}\label{eq:gb}
  \Vert g \Vert_2^2 \le C s  \Vert f\Vert_1, \; 
\Vert b \Vert_1 \le 2 \Vert f\Vert_1.
\end{equation}

\nd
Now let $f\in  L^1 \cap L^2$, then $ b= \sum b_k $ converges both a.e. and in
$L^1 \cap L^2$, by the definition of $b_k$ and properties (ii),
(iii), where 
$$
b_k(x) = \left\{
\begin{array}{ll}
f(x) -\frac{1}{|I_k|} \int_{I_k} f dy  & x\in I_k \\
0 &   x\notin I_k.
\end{array}\right.
$$
It follows from Lemma 6.4 and properties $(ii) $ and $(iv)$ that 
\begin{align}  
&\int_{\R\backslash D_s^*} |m(H) b(x) | dx \le \sum_k \int_{\R\backslash D_s^*} |m(H)
b_k(x) | \ud x \notag\\
 \le& \sum_k \int_{I_k} |b_k(y)| dy \int_{ \R\backslash I_k^*} |m(H)(x,y)-
m(H)(x,\overline{y}_k) | \ud x \label{eq:*bf}\\
\le& A \sum_k \int |b_k(y)| dy \le 2A  \Vert f\Vert_1, \notag
\end{align}  
where $D_s^* =  \cup_k I_k^* $,  with \\
$I_k^* = ( 
\overline{y}_k -2 t_k, \overline{y}_k +2 t_k) \cup (-\overline{y}_k -2 t_k, -\overline{y}_k +2 t_k) $.
\vspace{.11in}
Since $| D_s^*| \le 4 | D_s| $, from (\ref{eq:gb}) and (\ref{eq:*bf})  we obtain the weak $(1,1)$ estimate.
\hfill $\Box$

\vspace{.24in}
\nd
{\bf Proof of Theorem \ref{th:Besov-mult}.}  For $ g\in  L^2 \cap B_p^{\alpha , q} (H) $,
\begin{align*}
 \Vert m(H) g \Vert_{ B_p^{\alpha , q} }  = &  
\Vert \Phi(H) m(H) g
\Vert_p + \bigg\{ \sum_{j=1}^{\infty} (2^{j\alpha }  \Vert \vphi_j(H) m(H) g
\Vert_p )^q \bigg\}^{1/q} \\
 = &\Vert \{ 2^{j\alpha }  \vphi_j(H) m(H) g \} \Vert_{\ell^q
(L^p)}.
\end{align*}
Using 
$\vphi_j (H)=  \sum_{\nu =-1}^1 (\vphi\psi)_{j+\nu} (H) \vphi_j(H)$, with
convention $\phi_0=\Phi , \phi_{-1}= 0$,
 we have 
$$ \quad \Vert \{  2^{j\alpha }  \vphi_j(H) m(H) g \} \Vert_{\ell^q
(L^p)} \le C_{p,q} \sum_{\nu = -1}^1 \Big\{\sum_{j=0}^{\infty}  2^{j\alpha q} \Vert   m_{j+\nu}
(H) \vphi_j(H) g  \Vert_p^q \Big\}^{ 1/q} ,$$ 
where $m_j= m (\vphi\psi)_j$. Therefore it is sufficient to show that $m_j(H)
$ are uniformly bounded on $L^p, 1<p<\infty .$ But according to
Theorem 6.5 this is true because $m_j= m \psi_j$ verify the obvious
condition 
$$ |m_j^{(k)} (x) | \le C |x|^{-k} ,$$ 
$ k=0,1$, with $C$ independent of $j$. \hfill $\Box$ 
\vspace{0.25in}

\nd
{\bf Proof of Lemma \ref{lem:ker-size}.} Assuming $ \Vert z K_j (\cdot , y) \Vert_2 \le C \lam^{1/2}$,
Schwartz inequality gives
\begin{align*}  
&\int_{|z|>t} |K_j (x, y) | dx = \\
&\int_{ \{|x-y|>t \} \cap  \{
|x+y|>t \} } (\min_\pm |x\pm y |)^{-1}  \big\vert (\min_\pm |x
\pm y| ) K_j(x,y) \big\vert dx \\
  \le&\Big( \int_{ \{ |x-y>t\} \cap \{ |x+y|>t \}} (\min |x
\pm y|)^{-2} dx \Big)^{1/2} \Vert z K_j(\cdot , y) \Vert_2 
 \le  C t^{- \frac{1}{2}} \lam^{ \frac{1}{2} }.
\end{align*} 
Next we need to show $ \Vert z K_j (\cdot , y)\Vert_2 \le C \lam^{1/2}$.
Clearly, 
\begin{align*}
 &\Vert z m_j(H) (\cdot , y)\Vert_2 \le   \Vert z m_j(H) (x , y)
\chi_{ \{x>1\} }  \Vert_2 \\
+& \Vert z m_j(H) (x , y)\chi_{ \{x \le 1 \} }  \Vert_2 
 + \Vert z m_j(H) (x , y)\chi_{ \{x<-1\} }  \Vert_2.
\end{align*}
We can show that each of these three terms 
is $ \le C_{\eps} \lam^{1/2}$.
We shall prove the estimate for the first term only since 
the other two terms can be proved similarly.
The discussion is divided into three cases:  $\;y>1$, 
$  \;|y|\le 1$ and  $ \;y<-1$.
Again here we indicate the proof for the case $y>1$ only. 

Let $y>1$, $x>1$ and consider the high frequency case
$j>J:= 4+[2\log_2 \eps]$ first. 
Recall that $j>J\Leftrightarrow \lam^{-1}>4\eps$. 
\begin{align*}
 &2\pi \min_\pm |x\pm y| m_j(H) (x , y) = \min |x\pm y| \int_{ |\xi |
>2\eps } m_j(\xi^2) e(x,\xi)\bar{ e}(y,\xi)d\xi \\
=& z \int_{ \xi 
>2\eps }^+  m_j(\xi^2) |C_+|^2 e^{i (x-y)\xi} d\xi + z \int_{ \xi 
<-2\eps }^-  m_j(\xi^2) (e^{ix\xi} + C'_- e^{-ix \xi}) \overline{  e^{i
y\xi} + C'_- e^{-i y \xi} } \\
:=& I^+(x,y) + I^-(x,y).
\end{align*}
Integrating by parts we get 
\begin{align*}
&\vert  I^+ (x,y) \vert \le |x-y|\cdot \Big\vert \int^+  m_j(\xi^2) |C_+|^2
e^{i (x-y)\xi} \ud\xi \Big\vert \\
=&  \Big\vert \int^+  \frac{d}{d\xi} (m_j(\xi^2) |C_+|^2 )
e^{i (x-y)\xi} \ud \xi \Big\vert \\
=&\sqrt{2\pi} \left\vert { \left[\frac{d}{d\xi} (m_j(\xi^2) |C_+|^2 \chi_{\{\xi
>0\}} ) \right] }^\vee (x-y) \right\vert .
\end{align*}
By Plancherel formula, 
\begin{equation} \label{eq:m-I+}
\Vert   I^+ (x,y) \chi_{\{x 
>1 \}} \Vert_2  \le \sqrt{2\pi} \Vert \frac{d}{d\xi} (m_j(\xi^2) |C_+|^2 \chi_{\{\xi
>0\}} )   \Vert_2 \le C_{\eps} \lam^{1/2},
\end{equation} 
where we used if $ 1/2\lam \le |\xi| \le 1/\lam $,
$$\left\{
\begin{array}{ll}
m_j(\xi^2) &=O(1),\\
\frac{d}{d\xi} [m_j(\xi^2) ] &=O(1/\xi),\\
|C_+|^2 &= O(1),\\
\frac{d}{d\xi} \big( |C_+|^2\big) &=O(1/\xi^4).
\end{array}\right.
$$
Similarly, one can show that 
\begin{equation}\label{eq:m-I-}
\Vert   I^- (x,y) \chi_{\{x 
>1 \}} \Vert_2   \le C_{\eps} \lam^{1/2}.
\end{equation} 
Combing (\ref{eq:m-I+}), (\ref{eq:m-I-}), we get 
\begin{equation} 
 \Vert   z\, m_j(H)(x,y)  \chi_{\{x 
>1 \}} \Vert_2   \le C_{\eps} \lam^{1/2}.
\end{equation} 

Estimation for the low frequency case $ j\le J $ can be obtained by
following the same line (with a suitable modification when
necessary) for the high frequency case above,
 except that we use certain asymptotic properties near
the origin instead of $\infty$, (consult $\S 4$). 

We are left with the first inequality (\ref{eq:m-K_j}) concerning 
the ``size'' of
the kernel. The proof of (\ref{eq:m-K_j}) is similar to but 
easier than that of 
(\ref{eq:zker-size}) and may be left
as a dull exercise. This closes the proof of Lemma \ref{lem:ker-size}. \hfill $\Box$

\vspace{.22in}
\noindent
{\em Outline of the proof of Lemma \ref{lem:der-ker-size}}. 
  Lemma \ref{lem:der-ker-size} can be proved in the
same fashion as Lemma \ref{lem:ker-size}.  Assuming 
(\ref{eq:m-zder-K}), apply Schwartz 
inequality to get for $\textrm{all} \; y$ 
\begin{align*}
&\quad\int_{|z|>t } \left\vert  \frac{\partial}{\partial y} K_j(x,y) \right\vert dx 
\le\\ 
&\left( \int_{ \{|x-y|>t\}\cap \{
|x+y|>t\} } (\min |x\pm y|)^{-2} dx \right)^{1/2} \cdot \Vert z
\frac{\partial}{\partial y} K_j(\cdot, y) \Vert_2  
\le C t^{-1/2} \lam^{-1/2}. 
\end{align*}
Inequalities (\ref{eq:m-der-K}), (\ref{eq:m-zder-K}) measure the $L^2$ -norm of $
\frac{\partial}{\partial y}
K_j(\cdot ,y)$ and $ z \frac{\partial }{ \partial y}
K_j(\cdot ,y)$, which are derivative analogue of 
(\ref{eq:m-K_j}), (\ref{eq:zker-size}) in Lemma
\ref{lem:ker-size}. We only indicate here some points for 
(\ref{eq:m-zder-K}) since (\ref{eq:m-der-K}) is easier to deal with. 
Consider the high frequency case $j>J$ first. 
To
prove (\ref{eq:m-zder-K}) we break the function $x \mapsto z
\frac{\partial}{\partial y} K_j(x,y) $ into three parts: its
restriction to the sets $ \{ x>1 \} , \{ |x| \le 1 \}$ and  $\{x <-1
\}$. As before we are able to show that the $L^2$-norm of these
restrictions (in $x$) is  $\le C \lam^{-1/2}.$

For instance, in the case $y>1, x>1$, the identities     
$$
\left\{
\begin{array}{ll}
\frac{\partial}{\partial y} e_+ (y,\xi) & =i \xi  e_+ (y,\xi) \\
\frac{ \partial}{\partial y} e_- (y,\xi) & =i \xi ( e^{i y \xi}-{C'}_- e^{
-iy\xi} ) 
\end{array}\right.
$$
tell that the integral expression of $ z \frac{ \partial}{\partial y}
K_j(x,y) $ differs from that of $ z K_j(x,y) $ only by a factor $ i
\xi$ (up to a $\pm$ sign), for which reason we use the estimate $
\frac{d}{d\xi} [\xi m_j(\xi^2)] =O(1),\; \xi \rightarrow \infty (
1/2\lam \le |\xi| \le 1/\lam )$ in place of the estimate $ \frac{d}{d\xi} [ m_j(\xi^2)] =O(1/\xi).$ 

The interested reader can check the remaining cases as an exercise.
The corresponding inequality is valid for the low frequency case,
based on some simple asymptotic estimates as $ \xi\rightarrow 0.$
\hB

\section{Identification of $B^{\alpha, q}_p(H), 0 <\alpha <1$}

Generalized Besov space method was considered in \cite{J92}, \cite{JP85} and 
\cite{KRY97} in the study of perturbation of Schr\"odinger operators.
In application to PDE problems it is of interest to identify
these spaces.

 The spaces $B^{\alpha, q}_p(H)$ we have defined using ($\ref{besov-norm}$) and the system $\{\Phi, \vphi_j\}$ is  essentially of the same type as those  defined in \cite{JN} for $ p,q\ge 1$ and $\al\ge 0$. 
In \cite{JN}, sufficient conditions are given on $V$ so that $B_p^{\al,q}(H)$
can be identified with ordinary Besov spaces.
The proof is based on a real interpolation result, where the
interpolation spaces are defined via semigroups.
The following result is a variant of Theorem 5.1 in \cite{JN}.

Let $\cal{K}:=\{V:V=V_+ -V_-\;\textrm{so that}\; V_+ \in K_d^{loc},\; 
V_-\in K_d \}$, where $K_d$ denote the Kato class (see \S 1,\;\cite{JN} or 
\cite{Si82}).

\begin{theorem}\label{th:BZh}\quad Suppose $V\in\cal{K}$ and $\cal{D}(H^m)=W_p^{2m}$ for 
some $m\in \N $ and
$1\le p < \infty.$ Then for $1\le q \le \infty , 0<\al<m,
B^{\al,q}_p(H) =B^{2\al,q}_p(\R^d) $ (with equivalent norms).  
\end{theorem}

\vspace{.0508in}

Theorem 7.1 can be directly proved  
by following the proof of  Theorem 5.1 in \cite{JN} with obvious modifications.
Indeed, noting that $\mathcal{D}(H^m)=W_p^{2m}$, the 
proof is contained in the commutative diagram
$$
\begin{CD}
B^{\al,q}_p(H)@> = >> (L^p, \mathcal{D}(H^m) )_{\theta, q}, \\
@AAA         @AAA\\  
 B^{2\al,q}_p(\R^d)@> =>> (L^p, W_p^{2m} )_{\theta, q} \\
\end{CD}
$$  
with $ \theta =\frac{\al}{m}$. 

\vspace{.10in} 
\noindent
Remark 1. For convenience we state the theorem as above.
Note that the Besov norm was defined in \cite{JN} using $4$-adic system, 
while we have used dyadic system in this paper. 
By the way in
the condition $B(p,m)$ in \cite{JN}    $W_p^{m}$ should be  $W_p^{2m}$ as above.

\vspace{.110in}
\nd
Remark 2. Note that the condition of $V$ on the domain of $H^m$ is 
 equivalent
to Assumption $B(p, m)$ in \cite{JN}, which assumes that
 $(H+M)^{-m} $ is a bounded map from $L^p(\R^d)$ to $ W_p^{2m}(\R^d)$ with a
bounded inverse. 

It is essential to verify the domain condition on $H^m$
or, the assumption $B(p,m)$.  
In his communication to the second author A.~Jensen explained  that
it is easy to show that if $V$ is bounded
relative to $\triangle$ on $L^p (\R^d)$ with relative bound less than
one, then the condition $B(p, m)$ is satisfied for $m=1$.  
For $m>1$, the condition $B(p,m)$ is valid for all $m\ge
1$ and all $p$ 
if $V$ is $ C^\infty $ with all derivatives bounded.

In the following let $V$ be the barrier potential defined in $\S$1.
Obviously $V<<-\De$ with relative bound zero,
 satisfying the conditions in Theorem \ref{th:BZh}. Thus 
$B_p^{\alpha,q}(H)=B_p^{2\alpha,q}(\R) $ for $ 1\le p<\infty,
1\le q\leq \infty$, $0<\alpha< 1$. This, combined with 
Theorem \ref{th:Besov-mult} implies the following multiplier result
on ordinary Besov spaces. 

\begin{proposition}\label{prop:B(R)-mult} 
Suppose $m\in L^{\infty}$ be as in Theorem \ref{th:Besov-mult}. 
Then $m(H)$ is bounded on $B_p^{\alpha,q}(\R) $ for $ 1<p<\infty$,
$1\le q\leq \infty$, $0<\alpha< 2$. 
\end{proposition}

The other interesting result follows from the discussion above for barrier 
potential and Theorem 5.2 in \cite{JN}. Note that in one dimension we can take
equality for $\beta$.

\begin{proposition}\label{prop:Besov-e{-itH}} 
Suppose $ 1\le p<\infty$, $1\le q\leq \infty$ and $0<\alpha< 2-2\beta$ 
with $\beta =\vert \frac{1}{2}-\frac{1}{p}\vert$.
Then $e^{-itH}$ maps $B_p^{\alpha+2\beta,q}(\R)$ continuously
to $B_p^{\alpha,q}(\R) $.
Moreover, $e^{-itH}$ maps $B_p^{2\beta,q}(\R)$ continuously 
to $L^p$.
In both cases the operator norm is less than or equal to $C\langle t\rangle^\beta$.
\end{proposition}

For $m=2$, we have good reason in doubting the verification of 
the domain condition for $H^m$.

\vspace{.2in}
\noindent
{\bf Conjecture}.  $B^{\alpha, q}_p(H) \neq B^{\alpha,
q}_p(H_0)$, \quad  $ \alpha=2$.  

\vspace{.10in}
To see the reason we compare $H^2$ and $H_0^2$.
Write $H^2= H_0^2 +H_0V +VH_0 +V^2$.  If we take $A=H_0^2, B=H_0V +VH_0 +V^2$, 
the only term that could cause problem is the term $H_0V$, which
involves formally Dirac delta distributions and their first
derivatives. On  the other hand, \mbox{Theorem 3.2.2} in \cite{AK}
tells that the
domain of the operator $H_0 + c_1 \delta +c_2 \delta^\prime $ consists
of functions $u \in W_2^2(\R\setminus \{0\})$ with $u$ satisfying certain
boundary condition at the origin.   
So if $(H+ M)^2 $ is bounded from $W^{4,p}$ to $L^p,\; p=2 $, we would
have that the domain of the commutator $[V, H_0]$ is $W^{4,p}$, $p=2$, which
is not the case by the above theorem in \cite{AK}.

\end{document}